\documentclass{article}

\usepackage{amsfonts,amsmath,amssymb}
\usepackage{hyperref}
\usepackage{color}
\usepackage{graphicx}
\usepackage{subcaption}
\usepackage[section]{placeins}
\usepackage[utf8]{inputenc}
\usepackage{tikz}

\definecolor{darkblue}{rgb}{0,0.08,0.45}
\hypersetup{colorlinks,linkcolor=darkblue,citecolor=darkblue,filecolor=darkblue,urlcolor=darkblue}
\hypersetup{bookmarksopen,bookmarksnumbered}

\newcommand{\ipgi}[2]{\ensuremath{\left\langle #1 , #2 \right\rangle_{\surfaceint}}}
\newcommand{\norm}[1]{\ensuremath{\left\lVert #1 \right\rVert}}
\renewcommand{\d}{\ensuremath{\,\mathrm{d}}}
\newcommand{\dx}{\ensuremath{\d\mathbf{x}}}
\newcommand{\dy}{\ensuremath{\d\mathbf{y}}}
\newcommand{\normal}{\ensuremath{\hat{\mathbf{n}}}}
\newcommand{\exterior}{\ensuremath{\mathrm{ext}}}
\newcommand{\interior}{\ensuremath{\mathrm{int}}}
\newcommand{\ptot}{\ensuremath{p_{\mathrm{tot}}}}
\newcommand{\psca}{\ensuremath{p_{\mathrm{sca}}}}
\newcommand{\pinc}{\ensuremath{p_{\mathrm{inc}}}}
\newcommand{\pint}{\ensuremath{p_{\mathrm{int}}}}
\newcommand{\utot}{\ensuremath{u_{\mathrm{tot}}}}
\newcommand{\usca}{\ensuremath{u_{\mathrm{sca}}}}
\newcommand{\uinc}{\ensuremath{u_{\mathrm{inc}}}}
\newcommand{\uint}{\ensuremath{u_{\mathrm{int}}}}
\newcommand{\uext}{\ensuremath{u_{\mathrm{ext}}}}
\newcommand{\vint}{\ensuremath{\mathbf{v}_\interior}}
\newcommand{\vext}{\ensuremath{\mathbf{v}_\exterior}}
\newcommand{\wint}{\ensuremath{\mathbf{w}_\interior}}
\newcommand{\wext}{\ensuremath{\mathbf{w}_\exterior}}
\newcommand{\SLP}{\ensuremath{\mathcal{V}}}
\newcommand{\DLP}{\ensuremath{\mathcal{K}}}
\newcommand{\SLPint}{\ensuremath{\SLP_\interior}}
\newcommand{\DLPint}{\ensuremath{\DLP_\interior}}
\newcommand{\SLPext}{\ensuremath{\SLP_\exterior}}
\newcommand{\DLPext}{\ensuremath{\DLP_\exterior}}
\newcommand{\ID}{\ensuremath{I}}
\newcommand{\SL}{\ensuremath{V}}
\newcommand{\DL}{\ensuremath{K}}
\newcommand{\AD}{\ensuremath{T}}
\newcommand{\HS}{\ensuremath{D}}
\newcommand{\CA}{\ensuremath{A}}
\newcommand{\IDint}{\ensuremath{\ID_\interior}}
\newcommand{\SLint}{\ensuremath{\SL_\interior}}
\newcommand{\DLint}{\ensuremath{\DL_\interior}}
\newcommand{\ADint}{\ensuremath{\AD_\interior}}
\newcommand{\HSint}{\ensuremath{\HS_\interior}}
\newcommand{\CAint}{\ensuremath{\CA_\interior}}
\newcommand{\IDext}{\ensuremath{\ID_\exterior}}
\newcommand{\SLext}{\ensuremath{\SL_\exterior}}
\newcommand{\DLext}{\ensuremath{\DL_\exterior}}
\newcommand{\ADext}{\ensuremath{\AD_\exterior}}
\newcommand{\HSext}{\ensuremath{\HS_\exterior}}
\newcommand{\CAext}{\ensuremath{\CA_\exterior}}
\newcommand{\SLintm}{\ensuremath{\mathring{\SL}_\interior}}
\newcommand{\DLintm}{\ensuremath{\mathring{\DL}_\interior}}
\newcommand{\ADintm}{\ensuremath{\mathring{\AD}_\interior}}
\newcommand{\HSintm}{\ensuremath{\mathring{\HS}_\interior}}
\newcommand{\SLextm}{\ensuremath{\mathring{\SL}_\exterior}}
\newcommand{\DLextm}{\ensuremath{\mathring{\DL}_\exterior}}
\newcommand{\ADextm}{\ensuremath{\mathring{\AD}_\exterior}}
\newcommand{\HSextm}{\ensuremath{\mathring{\HS}_\exterior}}
\newcommand{\CAintscaled}{\ensuremath{\widehat{\CA}_\interior}}
\newcommand{\CAextscaled}{\ensuremath{\widehat{\CA}_\exterior}}
\newcommand{\MASS}{\ensuremath{\mathbf{M}}}
\newcommand{\MASSint}{\ensuremath{\MASS_\interior}}
\newcommand{\MASSext}{\ensuremath{\MASS_\exterior}}
\newcommand{\MASSintinv}{\ensuremath{\MASSint^{-1}}}
\newcommand{\MASSextinv}{\ensuremath{\MASSext^{-1}}}
\newcommand{\TR}{\ensuremath{Z}}
\newcommand{\TRint}{\ensuremath{\TR_\interior}}
\newcommand{\TRext}{\ensuremath{\TR_\exterior}}
\newcommand{\TRintinv}{\ensuremath{\TR_\interior^{-1}}}
\newcommand{\TRextinv}{\ensuremath{\TR_\exterior^{-1}}}
\newcommand{\TRextint}{\ensuremath{\TRintinv \TRext}}
\newcommand{\TRintext}{\ensuremath{\TRextinv \TRint}}
\newcommand{\PR}{\ensuremath{\mathbf{P}}}
\newcommand{\PRintext}{\ensuremath{\PR_\mathrm{int,ext}}}
\newcommand{\PRextint}{\ensuremath{\PR_\mathrm{ext,int}}}
\newcommand{\FEM}{\ensuremath{\mathcal{F}_\interior}}
\newcommand{\kint}{\ensuremath{k_\interior}}
\newcommand{\kext}{\ensuremath{k_\exterior}}
\newcommand{\rhoint}{\ensuremath{\rho_\interior}}
\newcommand{\rhoext}{\ensuremath{\rho_\exterior}}
\newcommand{\rhointext}{\ensuremath{\frac{\rhoint}{\rhoext}}}
\newcommand{\rhoextint}{\ensuremath{\frac{\rhoext}{\rhoint}}}
\newcommand{\cint}{\ensuremath{c_\interior}}
\newcommand{\cext}{\ensuremath{c_\exterior}}
\newcommand{\lambdaint}{\ensuremath{\lambda_\interior}}
\newcommand{\lambdaext}{\ensuremath{\lambda_\exterior}}
\newcommand{\hint}{\ensuremath{h_\interior}}
\newcommand{\hext}{\ensuremath{h_\exterior}}
\newcommand{\thetaint}{\ensuremath{\theta_\interior}}

\newcommand{\varphiint}{\ensuremath{\varphi_\interior}}
\newcommand{\varphiext}{\ensuremath{\varphi_\exterior}}
\newcommand{\phiint}{\ensuremath{\phi_\interior}}
\newcommand{\phiext}{\ensuremath{\phi_\exterior}}
\newcommand{\psiint}{\ensuremath{\psi_\interior}}
\newcommand{\psiext}{\ensuremath{\psi_\exterior}}
\newcommand{\Green}{\ensuremath{G}}
\newcommand{\Greenint}{\ensuremath{\Green_\interior}}
\newcommand{\Greenext}{\ensuremath{\Green_\exterior}}
\newcommand{\volume}{\ensuremath{\Omega}}
\newcommand{\volumeint}{\ensuremath{\volume_\interior}}
\newcommand{\volumeext}{\ensuremath{\volume_\exterior}}
\newcommand{\surface}{\ensuremath{\Gamma}}
\newcommand{\surfaceint}{\ensuremath{\surface_\interior}}
\newcommand{\surfaceext}{\ensuremath{\surface_\exterior}}
\newcommand{\Nint}{\ensuremath{N_\interior}}
\newcommand{\Next}{\ensuremath{N_\exterior}}
\newcommand{\traceD}{\ensuremath{\gamma_D}}
\newcommand{\traceN}{\ensuremath{\gamma_N}}
\newcommand{\traceDe}{\ensuremath{\traceD^+}}
\newcommand{\traceDi}{\ensuremath{\traceD^-}}
\newcommand{\traceNe}{\ensuremath{\traceN^+}}
\newcommand{\traceNi}{\ensuremath{\traceN^-}}

\title{The Boundary Element Method for Acoustic Transmission with Nonconforming Grids\footnote{\textcopyright~2024. This manuscript version is made available under the CC-BY-NC-ND 4.0
license. This manuscript is published in the Journal of Computational and Applied Mathematics in final
form at \url{https://doi.org/10.1016/j.cam.2024.115838}.}}
\author{Elwin van 't Wout\thanks{Institute for Mathematical and Computational Engineering, School of Engineering and Faculty of Mathematics, Pontificia Universidad Católica de Chile, Santiago, Chile. Contact: e.wout@uc.cl}}
\date{February 20, 2024}

\begin{document}

\maketitle

\begin{abstract}
	Acoustic wave propagation through a homogeneous material embedded in an unbounded medium can be formulated as a boundary integral equation and accurately solved with the boundary element method. The computational efficiency deteriorates at high frequencies due to the increase in mesh size with a fixed number of elements per wavelength and also at high material contrasts due to the ill-conditioning of the linear system. This study presents the design of boundary element methods feasible for nonconforming surface meshes at the material interface. The nonconforming algorithm allows for independent grid generation, improves flexibility, and reduces the degrees of freedom. It works for different boundary integral formulations for Helmholtz transmission problems, operator preconditioning, and coupling with finite element solvers. The extensive numerical benchmarks at canonical configurations and an acoustic foam model confirm the significant improvements in computational efficiency when employing the nonconforming grid coupling in the boundary element method.
\end{abstract}


\section{Introduction}

The boundary element method (BEM) is an efficient algorithm to numerically solve harmonic wave propagation in piecewise homogeneous materials. The algorithm allows for large-scale simulations in acoustics, electromagnetics, and elastodynamics, among other disciplines~\cite{lahaye2017modern, chew2001fast}. A rich literature on fractional Sobolev spaces supports its accuracy~\cite{nedelec2001acoustic, steinbach2008numerical, hsiao2008boundary, sauter2011boundary}. The key distinction between the BEM and numerical methods such as finite differences and the finite element method (FEM) is the reformulation of the volumetric differential equation into boundary integral equations at the material interfaces. Green's functions in the BEM guarantee accurate simulation of exterior scattering and wave transmission. Furthermore, accelerators such as the fast multipole method and hierarchical matrix compression efficiently compute large-scale problems on modern computer architectures. However, the BEM is limited to models that possess Green's functions, thus restricting the applicability primarily to piecewise homogeneous materials with a linear response to the wave field.

This study considers acoustic wave transmission through a homogeneous material embedded in free space. Such configurations can be modelled by boundary integral formulations of the Helmholtz equation for the exterior and interior regions. The transmission conditions at the material interface couple the field representations into a surface potential problem. The different choices of representation formulas and surface potentials lead to a wealth of boundary integral formulations with specific computational characteristics~\cite{wout2021benchmarking}. The boundary integral equations are then numerically solved with Galerkin or collocation methods. The Galerkin discretisation requires a surface mesh, typically made of triangular elements~\cite{smigaj2015solving}. Differently, collocation-based schemes such as Nyström discretisation~\cite{bruno2001fast, liu2019non}, isogeometric analysis~\cite{chen2022sample, chen2022bi}, and neural networks~\cite{anitescu2019artificial} do not require a triangularisation of the surface. This study considers Galerkin discretisation on triangular surface grids because of its meshing flexibility and mathematical solid foundation.

{\color{black}All boundary integral formulations for transmission problems involve a set of boundary integral operators based on the interior Green's function and others based on the exterior wavenumber. The interior and exterior boundary integral representations must be coupled by enforcing} some form of field continuity across the discrete representation of the material interface. This condition is relatively straightforward for conforming meshes, that is, a unique grid at each surface where the nodes, edges, and faces are consistent with the structure of neighbouring elements. {\color{black}However, this also means that the grid resolution has to be sufficiently fine to represent the smallest wavelength between the exterior and interior sides of each material interface.} The novelty of this manuscript is the Galerkin discretisation of the exterior and interior operators on independent triangular surface meshes, {\color{black}which are then coupled by nonconforming grid projections. The added flexibility to the BEM's design alleviates the mesh restrictions at the interface's side with longer wavelength} and improves the overall computational efficiency for transmission problems.

Nonconforming grids are widely used in numerical methods for partial differential equations such as the FEM. In those cases, two volumetric meshes for different subdomains are independently created so that the boundaries of these grids might not match at the common interface~\cite{flemisch2006elasto}. There are many stable and efficient algorithms to perform the data transfer between nonconforming meshes (cf.~\cite{farhat1998load, boer2007review}). Examples include interpolation~\cite{beckert2001multivariate}, projection on Lagrange multipliers~\cite{hansbo2005lagrange}, and mortar techniques~\cite{kloppel2011fluid}. Allowing for nonconforming grids is particularly beneficial to scenarios such as multiphysics models~\cite{keyes2013multiphysics}, contact problems~\cite{tezduyar2001finite}, domain decomposition~\cite{houzeaux2017domain}, and modular software design~\cite{degroote2013partitioned}. The nonconforming FEM also extends the feasible frequency range for acoustic transmission, compared to standard techniques with the same number of degrees of freedom~\cite{flemisch2012non}.

There is a sharp contrast between the extensive literature and the impressive developments of nonconforming domain decomposition methods for differential equations compared to boundary integral equations. Furthermore, the techniques developed for the FEM cannot be applied directly to the BEM due to the reduction of the volumetric wave field to a surface potential. First, most nonconforming BEM techniques in the literature consider {\color{black}a different approach than those adopted in this manuscript. In those studies, the surface mesh is partitioned into separate surface patches that may not match at one-dimensional contours. Such surface-based domain decomposition approaches include} mortar boundary elements~\cite{healey2010mortar, cools2012mortar, cools2015mortar}, discontinuous Galerkin discretisation~\cite{peng2015domain, han2017domain, chen2018nonconformal, kong2019discontinuous, mi2019convergence, liu2021massively, huang2022simplified}, Nitsche methods~\cite{chouly2012nitsche, dominguez2014posteriori, heuer2017nonconforming, dominguez2018posteriori}, overlapping partitions~\cite{dault2014generalized}, and multi-branch basis functions~\cite{huang2021multibranch}. {\color{black}Second, this study neither follows the multi-domain approach of partitioning volumetric subdomains~\cite{wu2008multi}, leading to additional transmission interfaces. In contrast, the proposed algorithm does not decompose any geometric structure at all but instead uses two different meshes at the same physical surface. The reason for creating two meshes at each surface is that different boundary integral operators at the same interface can be assembled on another mesh, depending on the wavenumber of the Green's function involved. They are then coupled together with nonconforming grid projections.} Techniques for such nonconforming BEM implementations have only recently been reported. A nonoverlapping domain decomposition method was designed for nonpenetrable objects~\cite{peng2012nonconformal} and extended to coatings~\cite{zhao2016solving}, transmission problems~\cite{jiang2014solving, zhao2017efie, zhao2017solving}, mixed basis functions~\cite{zhao2015fast}, multiple traces formulations~\cite{zhao2017fast, zhao2018multiple}, and numerical accelerators~\cite{zhao2015hybrid, guo2016ie}. The coupling matrices at the nonconforming interfaces are calculated with numerical quadrature on a union mesh or Lagrange multipliers~\cite{jiang2018flexible}. Notice that all these references to literature concern Maxwell's equations for electromagnetics, not the Helmholtz equation. For acoustics, tearing and interconnecting techniques~\cite{langer2003boundary} perform volumetric domain decomposition but were applied to conforming meshes only.

This study presents the novel application of nonconforming BEM to the Helmholtz model for acoustic wave propagation. The proposed methodology in this study is distinctive from its electromagnetic counterparts. For instance, different basis functions, different boundary integral formulations, and a different mortar technique at nonconforming meshes will be presented, along with extensive computational benchmarks. Furthermore, the nonconforming BEM for acoustic transmission will be extended to operator preconditioning for screen problems and FEM-BEM coupling for heterogeneous materials. In the latter case, nonconforming techniques exist for FEM-BEM algorithms~\cite{langer2005coupled, fritze2005fem, vouvakis2007domain, ruberg2008coupling, merz2009structural, peters2012structural, liang2019coupled, jia2019twofold} and the volume-surface integral equations~\cite{li2018vsie, zhu2019discontinuous, li2020solving}, but all references use different algorithms than those presented in this manuscript. The manuscript's scope is the algorithmic development of acoustic BEM on nonconforming grids, excluding the numerical analysis of the boundary integral formulations.

This study designs a BEM for acoustic transmission at piecewise homogeneous regions with nonconforming meshes at material interfaces. The boundary integral formulations will be presented in Section~\ref{sec:formulation} and the nonconforming coupling algorithm in Section~\ref{sec:mortarbem}. Section~\ref{sec:extensions} explains the extension to operator preconditioning and FEM-BEM coupling. Finally, the numerical benchmarks in Section~\ref{sec:results} showcase the computational benefits of the nonconforming BEM.

\section{Formulation}
\label{sec:formulation}

This study considers the propagation of harmonic acoustic waves in media with a linear response. A bounded domain~$\volumeint \subset \mathbb{R}^3$ with surface~$\surface$ is embedded in an unbounded exterior domain~$\volumeext \subset \mathbb{R}^3$. The boundary is Lipschitz continuous with outward pointing unit normal~$\normal$. The time dependency of wave propagation is extracted by assumming an $e^{-\imath\omega t}$ waveform, where $\omega$ denotes the angular frequency and $\imath$ the imaginary unit. The known incident wave field~$\pinc$ has a frequency~$f$. Let us denote the unknown acoustic field by $\ptot$ and define the scattered field as $\psca = \ptot - \pinc$. The mass density is denoted by $\rhoext$ and $\rhoint$, and the speed of sound by $\cext$ and $\cint$ in the exterior and interior domains, respectively. They are all assumed to be constant. The exterior and interior wavenumbers are given by $\kext = 2\pi f / \cext$ and $\kint = 2\pi f / \cint$, respectively. The acoustic pressure field in such a configuration can accurately be described by the Helmholtz system
\begin{equation} \label{eq:helmholtz}
	\begin{cases}
		\Delta \ptot + \kext^2 \ptot = 0, & \text{in } \volumeext; \\
		\Delta \ptot + \kint^2 \ptot = 0, & \text{in } \volumeint; \\
		\traceDe \ptot = \traceDi \ptot, & \text{at } \surface; \\
		\frac1\rhoext \traceNe \ptot = \frac1\rhoint \traceNi \ptot, & \text{at } \surface; \\
		\lim_{\mathbf{r} \to \infty} |\mathbf{r}| (\partial_{|\mathbf{r}|} \psca - \imath \kext \psca) = 0
	\end{cases}
\end{equation}
where $\Delta$ denotes the Laplace operator and
\begin{subequations}
	\begin{align}
		\traceDe \ptot(\mathbf{x}) &= \lim_{\mathbf{y} \to \mathbf{x}} \ptot(\mathbf{y}), && \mathbf{x} \in \surface, \quad \mathbf{y} \in \volumeext, \\
		\traceDi \ptot(\mathbf{x}) &= \lim_{\mathbf{y} \to \mathbf{x}} \ptot(\mathbf{y}), && \mathbf{x} \in \surface, \quad \mathbf{y} \in \volumeint, \\
		\traceNe \ptot(\mathbf{x}) &= \lim_{\mathbf{y} \to \mathbf{x}} \nabla\ptot(\mathbf{y}) \cdot \normal(\mathbf{x}), && \mathbf{x} \in \surface, \quad \mathbf{y} \in \volumeext, \\
		\traceNi \ptot(\mathbf{x}) &= \lim_{\mathbf{y} \to \mathbf{x}} \nabla\ptot(\mathbf{y}) \cdot \normal(\mathbf{x}), && \mathbf{x} \in \surface, \quad \mathbf{y} \in \volumeint,
	\end{align}
\end{subequations}
the Dirichlet and Neumann traces, respectively.

\subsection{Boundary integral operators}

Since the model considers a piecewise homogeneous material, the volumetric Helmholtz equation~\eqref{eq:helmholtz} can be rewritten into boundary integral equations on the material interface. More information on the design process of boundary integral formulations can be found in~\cite{wout2021benchmarking}, where the same notation is used. Details on the functional analysis can be found in textbooks such as~\cite{nedelec2001acoustic, steinbach2008numerical, hsiao2008boundary, sauter2011boundary}.

Since this study introduces nonconforming meshes at material interfaces, the interior and exterior sides of the object's surface need to be distinguished. For this purpose, let us define $\surfaceint$ and $\surfaceext$ the triangular surface meshes at $\surface$ that will be used for the interior and exterior boundary integral operators, respectively. Here, only polyhedral surfaces will be considered so that the physical surfaces spanned by the meshes coincide with the material interfaces.

The operators mapping from the surface to volume are given by
\begin{subequations}
	\begin{align}
		[\SLPint\psi](\mathbf{x}) &= \iint_{\surfaceint} \Greenint(\mathbf{x},\mathbf{y}) \psi(\mathbf{y}) \dy && \text{for } \mathbf{x} \in \volumeint; \\
		[\SLPext\psi](\mathbf{x}) &= \iint_{\surfaceext} \Greenext(\mathbf{x},\mathbf{y}) \psi(\mathbf{y}) \dy && \text{for } \mathbf{x} \in \volumeext; \\
		[\DLPint\phi](\mathbf{x}) &= \iint_{\surfaceint} \frac{\partial \Greenint(\mathbf{x},\mathbf{y})}{\partial\normal(\mathbf{y})} \phi(\mathbf{y}) \dy && \text{for } \mathbf{x} \in \volumeint; \\
		[\DLPext\phi](\mathbf{x}) &= \iint_{\surfaceext} \frac{\partial \Greenext(\mathbf{x},\mathbf{y})}{\partial\normal(\mathbf{y})} \phi(\mathbf{y}) \dy && \text{for } \mathbf{x} \in \volumeext;
	\end{align}
\end{subequations}
the single-layer and double-layer potential integral operators, respectively. Here,
\begin{subequations}
	\begin{align}
		\Greenint(\mathbf{x},\mathbf{y}) &= \frac{e^{\imath \kint |\mathbf{x}-\mathbf{y}|}}{4\pi|\mathbf{x}-\mathbf{y}|}, && \mathbf{x}, \mathbf{y} \in \volumeint, \; \mathbf{x} \ne \mathbf{y}; \\
		\Greenext(\mathbf{x},\mathbf{y}) &= \frac{e^{\imath \kext |\mathbf{x}-\mathbf{y}|}}{4\pi|\mathbf{x}-\mathbf{y}|}, && \mathbf{x}, \mathbf{y} \in \volumeext, \; \mathbf{x} \ne \mathbf{y}
	\end{align}
\end{subequations}
denote the Green's functions with the wavenumber of the interior and exterior region, respectively. Furthermore,
\begin{subequations}
	\begin{align}
		[\SLext \psi](\mathbf{x}) &= \iint_{\surfaceext} \Greenext (\mathbf{x},\mathbf{y}) \psi(\mathbf{y}) \d\mathbf{y} && \text{for } \mathbf{x} \in \surfaceext; \\
		[\DLext \phi](\mathbf{x}) &= \iint_{\surfaceext} \frac{\partial}{\partial \normal(\mathbf{y})} \Greenext (\mathbf{x},\mathbf{y}) \phi(\mathbf{y}) \d\mathbf{y} && \text{for } \mathbf{x} \in \surfaceext; \\
		[\ADext \psi](\mathbf{x}) &= \frac{\partial}{\partial \normal(\mathbf{x})} \iint_{\surfaceext} \Greenext (\mathbf{x},\mathbf{y}) \psi(\mathbf{y}) \d\mathbf{y} && \text{for } \mathbf{x} \in \surfaceext; \\
		[\HSext \phi](\mathbf{x}) &= -\frac{\partial}{\partial \normal(\mathbf{x})} \iint_{\surfaceext} \frac{\partial}{\partial \normal(\mathbf{y})} \Greenext (\mathbf{x},\mathbf{y}) \phi(\mathbf{y}) \d\mathbf{y} && \text{for } \mathbf{x} \in \surfaceext;
	\end{align}
\end{subequations}
which are called the exterior single-layer, double-layer, adjoint double-layer and hypersingular boundary integral operators, respectively. The interior boundary integral operators are defined equivalently but with the interior Green's function and the interior surface mesh. Remember that the normal is always pointing outwards, regardless of the operator.

With the distinction between interior and exterior meshes, the following identity mappings need to be introduced:
\begin{subequations}
	\label{eq:identityoperators}
	\begin{align}
		\IDext: \surfaceext &\to \surfaceext, \nonumber \\
		\varphiext &\mapsto \varphiext; \label{eq:idext} \\
		\IDint: \surfaceint &\to \surfaceint, \nonumber \\
		\varphiint &\mapsto \varphiint; \label{eq:idint} \\
		\TRext: \surfaceext &\to \surface, \nonumber \\
		\varphiext &\mapsto \varphiext; \label{eq:trext} \\
		\TRint: \surfaceint &\to \surface, \nonumber \\
		\varphiint &\mapsto \varphiint. \label{eq:trint}
	\end{align}
\end{subequations}
The first two operators are standard identity operators acting on the exterior or interior surface mesh, respectively. The last two operators are transmission operators that formally represent an identity mapping from the exterior or interior surface mesh towards the physical surface. Hence, the operator $\TRintinv \TRext$ maps a function on the exterior mesh towards itself on the interior mesh.

\subsection{Direct boundary integral formulations}

Let us consider the following definitions of the unknown interior and exterior fields:
\begin{subequations}
	\begin{align}
		\uint &= \begin{cases} 0 & \text{in } \volumeext, \\ \ptot & \text{in } \volumeint; \end{cases} \label{eq:uint} \\
		\uext &= \begin{cases} \ptot - \pinc & \text{in } \volumeext, \\ -\pinc & \text{in } \volumeint; \end{cases} \label{eq:uext}
	\end{align}
\end{subequations}
which are the interior field and the scattered field, respectively. These fields can be represented as
\begin{subequations}
	\label{eq:representation}
	\begin{align}
		\uint &= \SLPint \psiint - \DLPint \phiint & \text{in } \mathbb{R}^3 \setminus \surface, \label{eq:representation:int} \\
		\uext &= \SLPext \psiext - \DLPext \phiext & \text{in } \mathbb{R}^3 \setminus \surface, \label{eq:representation:ext}
	\end{align}
\end{subequations}
in terms of the unknown surface potentials
\begin{subequations}
	\label{eq:potentials:traces}
	\begin{align}
		\phiint &= \traceDi \ptot && \text{on } \surfaceint, \\
		\psiint &= \traceNi \ptot && \text{on } \surfaceint, \\
		\phiext &= -\traceDe \ptot && \text{on } \surfaceext, \\
		\psiext &= -\traceNe \ptot && \text{on } \surfaceext.
	\end{align}
\end{subequations}
Hence, the transmission conditions read
\begin{subequations}
	\label{eq:transmission:potentials}
	\begin{align}
		-\TRext \phiext &= \TRint \phiint && \text{on } \surface; \\
		-\frac1\rhoext \TRext \psiext &= \frac1\rhoint \TRint \psiint && \text{on } \surface.
	\end{align}
\end{subequations}
The traces of the representation formulas read
\begin{subequations}
	\label{eq:calderon:potentials}
	\begin{align}
		\left(\frac12\IDint - \CAint\right) \begin{bmatrix} \phiint \\ \psiint \end{bmatrix} &= \begin{bmatrix} 0 \\ 0 \end{bmatrix} && \text{on } \surfaceint, \label{eq:calderon:interior:potentials} \\
		\left(-\frac12\IDext - \CAext\right) \begin{bmatrix} \phiext \\ \psiext \end{bmatrix} &= \begin{bmatrix} \traceDe \uinc \\ \traceNe \uinc \end{bmatrix} && \text{on } \surfaceext, \label{eq:calderon:exterior:potentials}
	\end{align}
\end{subequations}
for the Calderón systems
\begin{subequations}
	\begin{align}
		\CAint &= \begin{bmatrix} -\DLint & \SLint \\ \HSint & \ADint \end{bmatrix}, \\
		\CAext &= \begin{bmatrix} -\DLext & \SLext \\ \HSext & \ADext \end{bmatrix}.
	\end{align}
\end{subequations}
Notice that the transmission conditions~\eqref{eq:transmission:potentials} yield two linear equations and the Calderón systems~\eqref{eq:calderon:potentials} another four linear equations, all with respect to four unknown potentials~\eqref{eq:potentials:traces}. The design of consistent boundary integral formulations follows different combinations of these model equations, as will be presented below.

\subsubsection{Multiple-traces formulation}

Substituting the identity operators in the Calderón systems~\eqref{eq:calderon:potentials} by the interface conditions~\eqref{eq:transmission:potentials} gives the system
\begin{equation}
	\label{eq:mtf}
	\begin{bmatrix}
		-\DLext & \SLext & \frac12 \TRextinv \TRint & 0 \\
		\HSext & \ADext & 0 & \frac12 \rhoextint \TRextinv \TRint \\
		-\frac12 \TRintinv \TRext & 0 & -\DLint & \SLint \\
		0 & -\frac12 \rhointext \TRintinv \TRext & \HSint & \ADint
	\end{bmatrix}
	\begin{bmatrix}
		\phiext \\ \psiext \\ \phiint \\ \psiint
	\end{bmatrix}
	= \begin{bmatrix}
		-\traceDe \uinc \\ -\traceNe \uinc \\ 0 \\ 0
	\end{bmatrix}
\end{equation}
which is called the multiple-traces formulation (MTF)~\cite{claeys2013multi}.

\subsubsection{Single-trace formulations}

The single-trace formulations use a single set of Dirichlet and Neumann traces at the boundary. Specifically,
\begin{subequations}
	\label{eq:singletrace:traces:ext}
	\begin{align}
		\phi &= -\phiext = \TRextinv \TRint \phiint, & \text{on } \surfaceext; \\
		\psi &= -\psiext = \rhoextint \TRextinv \TRint \psiint, & \text{on } \surfaceext,
	\end{align}
\end{subequations}
where the equalities follow from the transmission conditions. Substituting these traces in the interior Calderón system~\eqref{eq:calderon:interior:potentials} yields
\begin{align}
	\left(\frac12\IDint - \CAintscaled \right) \begin{bmatrix} \TRintinv \TRext \phi \\ \TRintinv \TRext \psi \end{bmatrix} &= \begin{bmatrix} 0 \\ 0 \end{bmatrix} \label{eq:calderon:stf:extint}
\end{align}
where
\begin{align}
	\CAintscaled &= \begin{bmatrix} -\DLint & \rhointext \SLint \\ \rhoextint \HSint & \ADint \end{bmatrix}
\end{align}
a scaled Calderón operator. Notice that the interior equation~\eqref{eq:calderon:stf:extint} maps potentials on the exterior mesh to the interior mesh. For consistency, and to accomodate combined field equations, let us consider
\begin{subequations}
	\label{eq:calderon:singletrace}
	\begin{align}
		\TRextinv \TRint \left(\frac12\IDint - \CAintscaled \right) \TRintinv \TRext \begin{bmatrix} \phi \\ \psi \end{bmatrix} &= \begin{bmatrix} 0 \\ 0 \end{bmatrix} && \text{on } \surfaceext, \\
		\left(\frac12\IDext + \CAext \right) \begin{bmatrix} \phi \\ \psi \end{bmatrix} &= \begin{bmatrix} \traceDe \uinc \\ \traceNe \uinc \end{bmatrix} && \text{on } \surfaceext,
	\end{align}
\end{subequations}
which are called the interior and exterior single-trace Calderón systems, respectively. They consist of four linear equations for two unknown potentials. Hence, linear combinations yield different boundary integral equations. Taking the difference of the exterior and interior Calderón systems yields
\begin{align*}
	\frac12 \left(\IDext - \TRextinv \TRint \IDint \TRintinv \TRext \right) \begin{bmatrix} \phi \\ \psi \end{bmatrix} + \left(\CAext + \TRextinv \TRint \CAintscaled \TRintinv \TRext \right) \begin{bmatrix} \phi \\ \psi \end{bmatrix} &= \begin{bmatrix} \traceDe \uinc \\ \traceNe \uinc \end{bmatrix}.
\end{align*}
The first term is zero (except for numerical errors) so that
\begin{align}
	\label{eq:pmchwt:ext}
	\left(\CAext + \TRextinv \TRint \CAintscaled \TRintinv \TRext \right) \begin{bmatrix} \phi \\ \psi \end{bmatrix} &= \begin{bmatrix} \traceDe \uinc \\ \traceNe \uinc \end{bmatrix} && \text{on } \surfaceext,
\end{align}
which is called the exterior PMCHWT formulation after the inventors of its electromagnetic variant: Poggio-Miller-Chang-Harrington-Wu-Tsai~\cite{poggio1973integral, chang1974surface, wu1977scattering-bor}, and also known as the Costabel-Stephan formulation~\cite{costabel1985direct}. Alternatively, taking the sum of the exterior and interior Calderón systems yields
\begin{align*}
	\frac12 \left(\IDext + \TRextinv \TRint \IDint \TRintinv \TRext \right) \begin{bmatrix} \phi \\ \psi \end{bmatrix} + \left(\CAext - \TRextinv \TRint \CAintscaled \TRintinv \TRext \right) \begin{bmatrix} \phi \\ \psi \end{bmatrix} &= \begin{bmatrix} \traceDe \uinc \\ \traceNe \uinc \end{bmatrix}.
\end{align*}
Assuming exact projection of the identity operator,
\begin{align}
	\label{eq:muller:ext}
	\IDext + \left(\CAext - \TRextinv \TRint \CAintscaled \TRintinv \TRext \right) \begin{bmatrix} \phi \\ \psi \end{bmatrix} &= \begin{bmatrix} \traceDe \uinc \\ \traceNe \uinc \end{bmatrix} && \text{on } \surfaceext,
\end{align}
which is called the exterior Müller formulation~\cite{muller1957grundprobleme}, and also known as the Kress-Roach formulation~\cite{kress1978transmission}.

\paragraph{Interior variant}
Other variants of single-trace formulations use the interior potentials
\begin{subequations}
	\label{eq:singletrace:traces:int}
	\begin{align}
		\phi &= \phiint = -\TRintinv \TRext \phiext, & \text{on } \surfaceint; \\
		\psi &= \psiint = -\rhointext \TRintinv \TRext \psiext, & \text{on } \surfaceint.
	\end{align}
\end{subequations}
Then, the interior PMCHWT formulation reads
\begin{align}
	\label{eq:pmchwt:int}
	\left(\TRintinv \TRext \CAextscaled \TRextinv \TRint + \CAint \right) \begin{bmatrix} \phi \\ \psi \end{bmatrix} &= \TRintinv \TRext \begin{bmatrix} \traceDe \uinc \\ \rhointext \traceNe \uinc \end{bmatrix} && \text{on } \surfaceint,
\end{align}
and
\begin{align}
	\label{eq:muller:int}
	\IDint + \left(\TRintinv \TRext \CAextscaled \TRextinv \TRint - \CAint \right) \begin{bmatrix} \phi \\ \psi \end{bmatrix} &= \TRintinv \TRext \begin{bmatrix} \traceDe \uinc \\ \rhointext \traceNe \uinc \end{bmatrix} && \text{on } \surfaceint,
\end{align}
is the interior Müller formulation, where
\begin{align*}
	\CAextscaled &= \begin{bmatrix} -\DLext & \rhoextint \SLext \\ \rhointext \HSext & \ADext \end{bmatrix}
\end{align*}
a scaled Calderón operator.

\subsection{Indirect boundary integral formulations}

The representation formulas~\eqref{eq:representation} directly give boundary integral formulations in terms of unknown surface potentials that represent acoustic field traces. Indirect representation formulas follow a different approach and define surface potentials that do not necessarily have a direct physical interpretation, here denoted by $\mu$ and $\nu$. Among the many options (cf.~\cite{wout2021benchmarking}), let us consider the high-contrast formulations~\cite{wout2022highcontrast}. In the following, superscript circles denote passage through the opposite mesh as
\begin{subequations}
	\begin{align}
		\mathring{X}_\exterior &= \TRextint X_\exterior \TRintext, \\
		\mathring{X}_\interior &= \TRintext X_\interior \TRextint
	\end{align}
\end{subequations}
for $X$ any of the boundary integral operators. Now, the representation formulas
\begin{subequations}
	\begin{align}
		\uint &= \SLPint \psiint - \DLPint \phiint, \\
		\uext &= \SLPext \mu_\exterior, \label{eq:representation:indirect:ext:slp}
	\end{align}
\end{subequations}
yield the system
\begin{align}
	\begin{bmatrix} \frac12\IDext - \ADintm & -\rhoextint \HSintm \SLext \\ \IDext & \frac12\IDext - \ADext \end{bmatrix}
	\begin{bmatrix} \traceNe \ptot \\ \mu_\exterior \end{bmatrix}
	&= \begin{bmatrix} \rhoextint \HSintm \traceDe \uinc \\ \traceNe \uinc \end{bmatrix} && \text{on } \surfaceext, \label{eq:highcontrast:ext:neu}
\end{align}
called the exterior high-contrast Neumann formulation. The representation formulas
\begin{subequations}
	\begin{align}
		\uint &= \SLPint \psiint - \DLPint \phiint, \\
		\uext &= -\DLPext \nu_\exterior \label{eq:representation:indirect:ext:dlp}
	\end{align}
\end{subequations}
yield the system
\begin{align}
	\begin{bmatrix} \frac12\IDext + \DLext & \IDext \\ -\rhointext \SLintm \HSext & \frac12\IDext + \DLintm \end{bmatrix}
	\begin{bmatrix} \nu_\exterior \\ \traceDe \utot \end{bmatrix}
	&= \begin{bmatrix} \traceDe \uinc \\ \rhointext \SLintm \traceNe \uinc \end{bmatrix} && \text{on } \surfaceext, \label{eq:highcontrast:ext:dir}
\end{align}
called the exterior high-contrast Dirichlet formulation. The representation formulas
\begin{subequations}
	\begin{align}
		\uint &= \SLPint \mu_\interior, \label{eq:representation:indirect:int:slp} \\
		\uext &= \SLPext \psiext - \DLPext \phiext
	\end{align}
\end{subequations}
yield the system
\begin{align}
	&\begin{bmatrix} \frac12\IDint + \ADextm & \rhointext \HSextm \SLint \\ -\IDint & \frac12\IDint + \ADint \end{bmatrix}
	\begin{bmatrix} \traceNi \utot \\ \mu_\interior \end{bmatrix} \nonumber \\
	&= \begin{bmatrix} \rhointext \TRextint \left(\HSext \traceDe \uinc + \left(\frac12 \IDext + \ADext\right) \traceNe \uinc\right) \\ 0 \end{bmatrix} && \text{on } \surfaceint,
	\label{eq:highcontrast:int:neu}
\end{align}
called the interior high-contrast Neumann formulation. Finally, the representation formulas
\begin{subequations}
	\begin{align}
		\uint &= -\DLPint \nu_\interior \label{eq:representation:indirect:int:dlp}, \\
		\uext &= \SLPext \psiext - \DLPext \phiext
	\end{align}
\end{subequations}
yield the system
\begin{align}
	&\begin{bmatrix} \frac12\IDint - \DLint & -\IDint \\ \rhoextint \SLextm \HSint & \frac12\IDint - \DLextm \end{bmatrix}
	\begin{bmatrix} \nu_\interior \\ \traceDi \utot \end{bmatrix} \nonumber \\
	&= \begin{bmatrix} 0 \\ \TRextint \left(\SLext \traceNe \uinc + \left(\frac12\IDext - \DLext\right) \traceDe \uinc\right) \end{bmatrix} && \text{on } \surfaceint,
	\label{eq:highcontrast:int:dir}
\end{align}
called the interior high-contrast Dirichlet formulation.

\subsection{Numerical discretisation}

The discretisation of the boundary integral operators follows a standard Galerkin method~\cite{smigaj2015solving}. All test and basis functions are continuous piecewise linear (P1) {\color{black}Lagrange polynomials}, with value one in a specific grid node and zero in all other nodes of the triangular surface mesh.

\section{Nonconforming grid projections}
\label{sec:mortarbem}

\begin{figure}[!ht]
	\centering
	\includegraphics[width=.7\columnwidth]{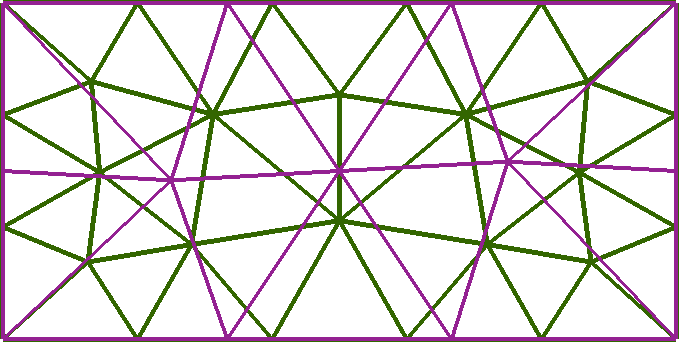}
	\caption{An example of two nonconforming triangular grids generated on the same rectangular domain.}
	\label{fig:mesh}
\end{figure}

This study considers the acoustic transmission at a homogeneous domain with surface~$\surface$, where two surface meshes $\surfaceint$ and $\surfaceext$ are independently generated and the interior and exterior boundary integral operators are assembled on the respective triangulations. Each individual mesh is assumed to be conforming and to cover the physical surface~$\surface$ exactly. However, the grid $\surfaceint$ might not be conforming with $\surfaceext$, as depicted in Figure~\ref{fig:mesh}. Hence, let us denote the number of nodes in the interior and exterior mesh by $\Nint$ and $\Next$, respectively. They are also the number of degrees of freedom in the discrete P1 spaces.

\subsection{Mortar matrices}

Since all dense boundary integral operators in the models are defined with test and basis functions on the same mesh, these can be assembled with a standard weak formulation. Differently, the transmission operators~\eqref{eq:trext} and~\eqref{eq:trint} map between different meshes on the same surface. For example, one formally has
\begin{align}
	[\TRint\varphiint](\mathbf{x}) = [\TRext\varphiext](\mathbf{x}) \text{ for } \mathbf{x} \in \surface,
\end{align}
where $\varphiint$ and $\varphiext$ are functions with the same values but defined on the meshes $\surfaceint$ and $\surfaceext$, respectively. Within the boundary integral formulations, these mortar operators are present in the form of $\varphiint = \TRextint \varphiext$, which represents an identity mapping from the exterior to the interior mesh. The weak formulation of this equation reads
\begin{align}
	\ipgi{\IDint\varphiint}{\thetaint} = \ipgi{\TRextint \varphiext}{\thetaint}
\end{align}
with a standard $L^2$ inner product. The interior identity operator was included for clarity. The discrete version, i.e., substituting basis and test functions, of this weak formulation is a set of linear equations given by
\begin{align}
	\MASSint \vint = \PRintext \vext
\end{align}
where $\vint$ and $\vext$ denote coefficient vectors on the interior and exterior mesh, respectively. The interior mass matrix $\MASSint$ has size $\Nint \times \Nint$ and element $(i,j)$ corresponds to the inner product of test function~$i$ and basis function~$j$ on the interior mesh. Differently, the mortar matrix $\PRintext$ has size $\Nint \times \Next$ and element $(i,j)$ corresponds to the inner product of test function~$i$ on the interior mesh and basis function~$j$ on the exterior mesh. These inner products are well defined since both meshes share the same physical surface. Finally, $\PRextint$ denotes the mortar matrix that maps from the interior to the exterior mesh and is the transpose of $\PRintext$.

\subsection{Advancing front method}

Volumetric domain decomposition techniques often require the computation of the same mortar matrix $\PRintext$. As mentioned in the introduction, many different algorithms exist to compute these matrices. Here, the \emph{Projection Algorithm for Nonmatching Grids} (PANG) will be used~\cite{gander2013algorithm}. This is an advancing front method that calculates the mortar matrix of P1 elements on nonconforming triangular surface meshes. This algorithm advances through the nonconforming grids by iteratively calculating triangle intersections, mortar contributions, and candidate neighbours. The reasons to use this specific algorithm are as follows. Firstly, open-source code of the algorithm is available~\cite{gander2013algorithm}. Secondly, the algorithm has linear computational complexity~\cite{gander2009algorithm}. Thirdly, it is robust in finite precision~\cite{mccoid2022provably}. Fourthly, the necessary triangle connectivity tables are already available in the standard mesh formats used in the BEM. Here, no modifications to the PANG were made since it directly calculates the mortar matrices~$\PRintext$ and~$\PRextint$.

\subsection{Discrete boundary integral formulations}

With the mortar matrices available, the discretised boundary integral formulations can readily be obtained. The discrete MTF~\eqref{eq:mtf} reads
\begin{equation}
	\label{eq:mtf:discrete}
	\begin{bmatrix}
		-\mathbf\DLext & \mathbf\SLext & \frac12 \PRextint & 0 \\
		\mathbf\HSext & \mathbf\ADext & 0 & \frac12 \rhoextint \PRextint \\
		-\frac12 \PRintext & 0 & -\mathbf\DLint & \mathbf\SLint \\
		0 & -\frac12 \rhointext \PRintext & \mathbf\HSint & \mathbf\ADint
	\end{bmatrix}
	\begin{bmatrix}
		\vext \\ \wext \\ \vint \\ \wint
	\end{bmatrix}
	= \begin{bmatrix}
		-\mathbf{f} \\ -\mathbf{g} \\ \mathbf{0} \\ \mathbf{0}
	\end{bmatrix}
\end{equation}
where the boldface operators denote the weak formulation of the corresponding boundary integral operators. The discrete PMCHWT formulation~\eqref{eq:pmchwt:ext} reads
\begin{equation}
	\label{eq:pmchwt:discrete}
	\begin{bmatrix}
		-\mathbf\DLext - \mathbf{\widehat{\DLint}} & \mathbf\SLext + \mathbf{\widehat{\SLint}} \\
		\mathbf\HSext + \mathbf{\widehat{\HSint}} & \mathbf\ADext + \mathbf{\widehat{\ADint}} \\
	\end{bmatrix}		
	\begin{bmatrix} \vext \\ \wext \end{bmatrix}
	= \begin{bmatrix} \mathbf{f} \\ \mathbf{g} \end{bmatrix}
\end{equation}
where the superscript hat denotes the transformation
\begin{align}
	\mathbf{\widehat{\SLint}} = \mathbf{\PRextint \MASSintinv \SLint \MASSintinv \PRintext}.
\end{align}
Notice that the inverse mass matrices are necessary for coupling the discrete spaces of the operator products (cf.~\cite{betcke2020product}). The discrete high-contrast formulation~\eqref{eq:highcontrast:ext:neu} reads
\begin{align}
	\begin{bmatrix} \frac12\MASSext - \mathbf{\widehat{\ADint}} & -\rhoextint \mathbf{\widehat{\HSint}} \MASSextinv \mathbf\SLext \\ \MASSext & \frac12\MASSext - \mathbf\ADext \end{bmatrix}
	\begin{bmatrix} \traceNe \ptot \\ \psiext \end{bmatrix}
	&= \begin{bmatrix} \rhoextint \mathbf{\widehat{\HSint}} \mathbf{f} \\ \mathbf{g} \end{bmatrix}.
\end{align}
All other formulations can be discretised similarly.

The discrete formulations mentioned above are all weak forms. Mass-matrix preconditioning yields the strong form, i.e., multiplying the discrete equations from the left with the matrices
\begin{equation*}
	\begin{bmatrix}
		\MASSextinv & 0 \\
		0 & \MASSextinv
	\end{bmatrix}
\end{equation*}
for the single-trace and high-contrast formulations and
\begin{equation*}
	\begin{bmatrix}
		\MASSextinv & 0 & 0 & 0 \\
		0 & \MASSextinv & 0 & 0 \\
		0 & 0 & \MASSintinv & 0 \\
		0 & 0 & 0 & \MASSintinv
	\end{bmatrix}
\end{equation*}
for the MTF. Alternatively, more elaborate preconditioning strategies such as Calderón and OSRC preconditioning~\cite{wout2021benchmarking} can be applied as usual.

\section{Extensions}
\label{sec:extensions}

Section~\ref{sec:formulation} presented the nonconforming BEM for acoustic transmission at a single homogeneous material embedded in free space. The extension to multiple scattering at disjoint objects is straightforward. Another use case of nonconforming BEM is efficient parameter studies: when the wavespeed in one of the subdomains changes, one only needs to reassemble the operators in that domain while keeping all other boundary integral operators. This section explains other interesting applications: FEM-BEM coupling and operator preconditioning.

\subsection{FEM-BEM coupling}
\label{sec:fembem}

Let us consider a heterogeneous material embedded in an unbounded homogeneous medium. The Helmholtz equation in the interior is given by
\begin{equation}
	\label{eq:helmholtz:heterogeneous}
	-\rhoint \nabla \cdot \left( \frac1\rhoint \nabla \pint \right) - \kint^2 \pint = 0
\end{equation}
where $\rhoint(\mathbf{x})$ and $\kint(\mathbf{x}) = 2\pi f / \cint(\mathbf{x})$ for $\mathbf{x} \in \volumeint$ smooth functions. Its weak formulation (cf.~\cite{steinbach2008numerical}) reads
\begin{align}
	&\iiint_{\volumeint} \frac1{\rhoint(\mathbf{x})} \nabla \pint(\mathbf{x}) \nabla \left(\rhoint(\mathbf{x}) q(\mathbf{x})\right) \dx - \iint_{\surfaceint} \traceNi\pint(\mathbf{x}) \traceDi q(\mathbf{x}) \dx \nonumber \\ &\qquad - \iiint_{\volumeint} \kint^2(\mathbf{x}) \pint(\mathbf{x}) q(\mathbf{x}) \dx = 0,
	\label{eq:weak:volume}
\end{align}
where $q$ denotes the test function. The discretisation uses node-based continuous piecewise linear (P1) Lagrange basis and test functions on a tetrahedral grid. The exterior Calderón system~\eqref{eq:calderon:exterior:potentials} yields
\begin{equation}
	\left(\frac12\IDext - \DLext\right) \traceDe\ptot + \SLext \traceNe\ptot = \traceDe\pinc.
\end{equation}
Defining the trace variable $\vartheta = \traceNe\ptot$ and using the transmission conditions~\eqref{eq:helmholtz}, the continuous formulations read
\begin{equation}
	\label{eq:fembem:standard}
	\begin{bmatrix}
		\FEM & -\frac\rhoint\rhoext \IDint \\
		\tfrac12\IDext - \DLext & \SLext \end{bmatrix}
	\begin{bmatrix} \pint \\ \vartheta \end{bmatrix}
	= \begin{bmatrix}
		0 \\
		\traceDe\pinc
	\end{bmatrix}
\end{equation}
where $\FEM$ denotes the volumetric part of the the weak formulation~\eqref{eq:weak:volume} of the Helmholtz equation inside the heterogeneous subdomain. This FEM-BEM coupling is known as the Johnson-Nédélec variant~\cite{johnson1980coupling, sayas2009validity}. The discrete formulation on nonconforming meshes reads
\begin{equation}
	\label{eq:fembem}
	\begin{bmatrix}
		\FEM & -\frac\rhoint\rhoext F_\mathrm{e} \TRextint \\
		\left(\tfrac12\IDext - \DLext\right) \TRintext F_\mathrm{r} & \SLext \end{bmatrix}
	\begin{bmatrix} \pint \\ \vartheta \end{bmatrix}
	= \begin{bmatrix}
		0 \\
		\traceDe\pinc
	\end{bmatrix}
\end{equation}
where the operators $\TRext$ and $\TRint$ are the mortar operators~\eqref{eq:identityoperators} between the exterior mesh and the boundary of the volumetric mesh. Furthermore, the operators $F_\mathrm{e}$ and $F_\mathrm{r}$ denote the extension and restriction from the degrees of freedom in the interior to the degrees of freedom on the boundary of the volumetric mesh. The nonconforming technique can easily be extended to other formulations, including stabilised ones~\cite{wout2022fembem, casenave2014coupled}.

\subsection{Operator preconditioning}
\label{sec:opprec}

Operator preconditioning (cf.~\cite{kirby2010functional, hiptmair2006operator}) effectively improves the convergence of linear solvers for the BEM~\cite{wout2021benchmarking, search2022towards} and the boundary integral operators for the preconditioner are independently assembled from the model. For instance, Calderón and OSRC preconditioners for acoustic transmission~\cite{wout2022pmchwt} can be assembled on a coarser mesh than the boundary integral operators for the model itself. This idea is not restricted to acoustic transmission and is also valid for impenetrable objects discretised with strategies such as the OSRC-preconditioned Burton-Miller formulation~\cite{wout2015fast} and opposite-order preconditioning of single-layer and hypersingular operators~\cite{steinbach1998construction}.

As an example, let us consider a Neumann screen problem
\begin{equation}
	\label{eq:neumannscreen}
	\begin{cases}
		\Delta \psca + \kext^2 \psca = 0, & \text{in } \volume; \\
		\traceN \psca = -\traceN \pinc, & \text{on } \surface; \\
		\lim_{\mathbf{r} \to \infty} |\mathbf{r}| (\partial_{|\mathbf{r}|} \psca - \imath \kext \psca) = 0;
	\end{cases}
\end{equation}
where~$\surface$ is an open two-dimensional manifold in the three-dimensional domain~$\volume$. The boundary integral formulation reads
\begin{equation}
	\begin{cases}
		\HSext \phi = \traceN\uinc, & \text{on } \surface; \\
		\usca = \DLPext \phi, & \text{in } \volume.
	\end{cases}
\end{equation}
Opposite-order preconditioning suggests the use of the single-layer integral operator as preconditioner~\cite{steinbach1998construction}, even in the case of screen problems where boundary singularities occur~\cite{bruno2013high}. The preconditioner does not need to have the same wavenumber or be assembled with the same numerical parameters as the original matrix~\cite{boubendir2015integral, escapil2021bi}. Hence,
\begin{equation}
	\TRintext \SLint \TRextint \HSext \phi = \TRintext \SLint \TRextint \traceN\uinc
\end{equation}
is a valid preconditioning strategy, where the model is assembled on an `exterior' mesh and the preconditioner on a coarser `interior' mesh.

\section{Results}
\label{sec:results}

This section presents numerical benchmarks of the nonconforming BEM on canonical test cases and a large-scale model of an acoustic foam.

\subsection{Computational framework}

The weak formulation of the BEM uses local P1 elements on a triangular surface mesh. Even though accelerators can readily be applied to the proposed technology, dense matrix arithmetic is used to limit the parameter space of the benchmarks. {\color{black}The surface integrals are numerically evaluated with quadrature rules of order four, with semi-analytical transformations to handle singularities~\cite{sauter2011boundary}.} The GMRES algorithm~\cite{saad1986gmres} solves the linear system, without restart and a termination criterion of $10^{-5}$. Table~\ref{table:formulations} summarises the boundary integral formulations for this study and mass matrix preconditioning is always applied since it improves the convergence of GMRES with little computational overhead.

The BEM was implemented with Bempp-cl (version~0.2.4)~\cite{betcke2021bempp} and the FEM with Fenics (version 2019.1)~\cite{logg2012automated}. The meshes were generated with Gmsh (version~4.9.5)~\cite{geuzaine2009gmsh}. The open-source Matlab code of the PANG algorithm~\cite{gander2013algorithm} and the PANG2 patch~\cite{mccoid2022provably} were rewritten into Python for compatibility with Bempp-cl and Fenics.

Shared-memory parallelisation was automatically performed through the PyOpenCL (version~2022.1)~\cite{numba} implementation of Bempp-cl for matrix assembly, Numpy (version~1.21.5)~\cite{numpy} for dense matrix algebra, Scipy (version~1.8)~\cite{scipy} for sparse matrix algebra, Numba (version~0.55.1)~\cite{numba} for the PANG, and Joblib (version~1.1)~\cite{joblib} for independent mortar matrices. All simulations were performed on a computer with two Intel(R) Xeon(R) CPU E5-2640~v4 @ 2.40~GHz processors, 20~cores, 40~threads, and 512~GB shared memory.

\begin{table}[!ht]
	\caption{Overview of the boundary integral formulations {\color{black}for acoustic transmission problems. The version refers to the unknown surface potentials being exterior or interior traces}. The number of operators refers to dense boundary integral operators only. The number of potentials refers to the unknown surface potentials.}
	\label{table:formulations}
	\centering
	\begin{tabular}{llrrr}
		\hline\hline
		formulation & version & equation & \#operators & \#potentials \\
		\hline
		PMCHWT & exterior & \eqref{eq:pmchwt:ext} & 8 & 2 \\
		PMCHWT & interior & \eqref{eq:pmchwt:int} & 8 & 2 \\
		Müller & exterior & \eqref{eq:muller:ext} & 8 & 2 \\
		Müller & interior & \eqref{eq:muller:int} & 8 & 2 \\
		multiple-traces & - & \eqref{eq:mtf} & 8 & 4 \\
		high-contrast & exterior Neumann & \eqref{eq:highcontrast:ext:neu} & 4 & 2 \\
		high-contrast & exterior Dirichlet & \eqref{eq:highcontrast:ext:dir} & 4 & 2 \\
		high-contrast & interior Neumann & \eqref{eq:highcontrast:int:neu} & 4 & 2 \\
		high-contrast & interior Dirichlet & \eqref{eq:highcontrast:int:dir} & 4 & 2 \\
		\hline\hline
	\end{tabular}
\end{table}

\subsection{Projection accuracy}
\label{sec:results:mortarerror}

The mortar matrices introduce a projection error when mapping between functions defined on P1 elements at nonconforming grids. On a continuous level, the projection operators~\eqref{eq:identityoperators} satisfy
\begin{align*}
	\TRextint \TRintext &= \IDint, \\
	\TRintext \TRextint &= \IDext.
\end{align*}
Hence, let us define the projection error of the mortar algorithm as
\begin{subequations}
	\label{eq:error:projection}
	\begin{align}
		E_\interior &= \norm{\mathbf\IDint - \MASSintinv \PRintext \MASSextinv \PRextint}, \label{eq:error:projection:int} \\
		E_\exterior &= \norm{\mathbf\IDext - \MASSextinv \PRextint \MASSintinv \PRintext}, \label{eq:error:projection:ext}
	\end{align}
\end{subequations}
corresponding to the interior and exterior mesh, respectively.

When both meshes are exactly the same, the projection error should be zero, except for rounding errors. Therefore, when a nonconforming mesh converges to the other mesh, the projection error is expected to converge to zero. Let us test this hypothesis on a triangular surface mesh for the unit square $[0,1] \times [0,1]$. The `interior' grid is fixed and `exterior' grids are generated from random perturbations of the nodes, drawn from a $\mathcal{N}(0,\sigma)$ distribution. To keep the same physical geometry, the inner nodes are perturbed in two dimensions, the boundary nodes only along the boundary, and the corner nodes not at all.

\begin{figure}[!ht]
	\centering
	\begin{subfigure}[b]{\columnwidth}
		\includegraphics[width=\columnwidth]{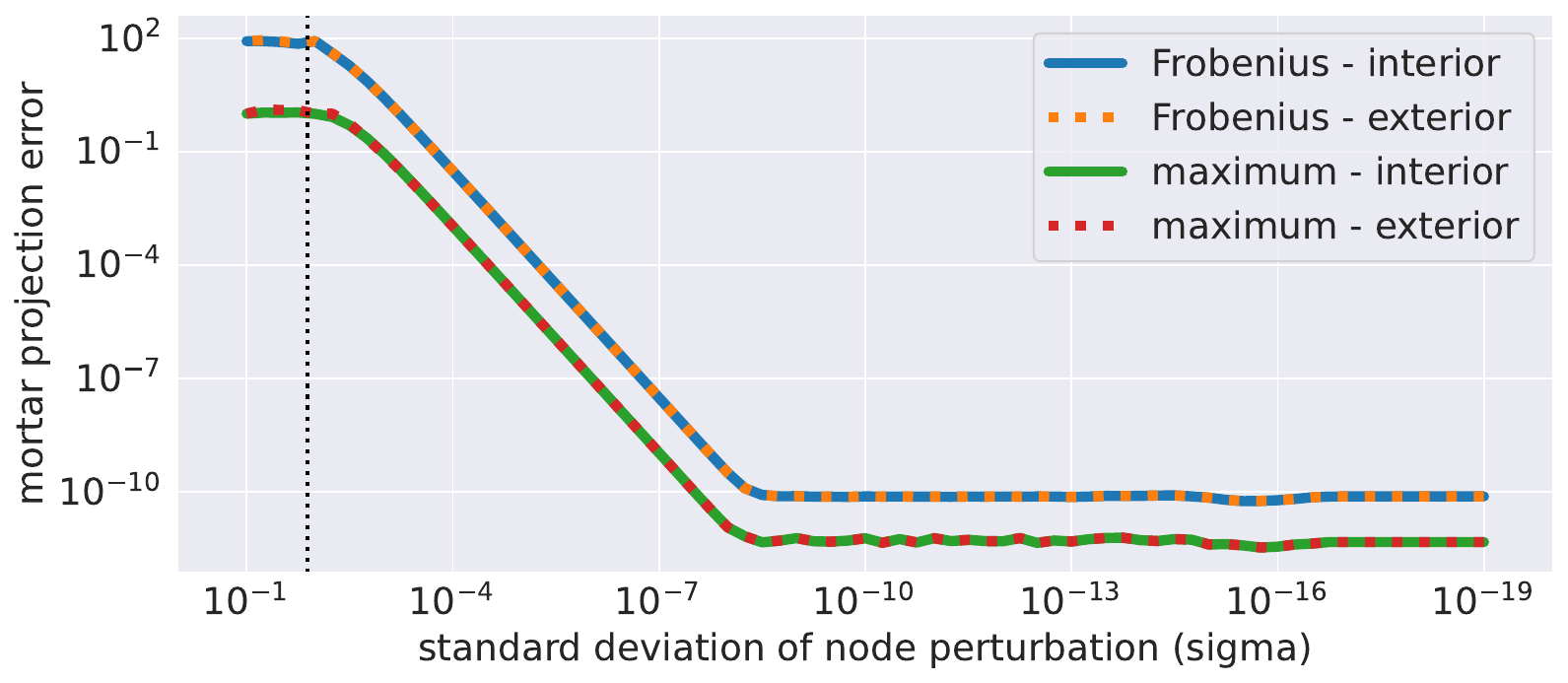}
		\caption{The projection error~\eqref{eq:error:projection} with respect to node perturbation with Gaussian noise having standard deviation $\sigma$ in the `exterior' mesh. The vertical line visualises the mesh width of $h=0.013$ of the fixed `interior' mesh, corresponding to 7044 nodes.}
		\label{fig:mortarerror:perturbations}
	\end{subfigure}
	\begin{subfigure}[b]{\columnwidth}
		\includegraphics[width=\columnwidth]{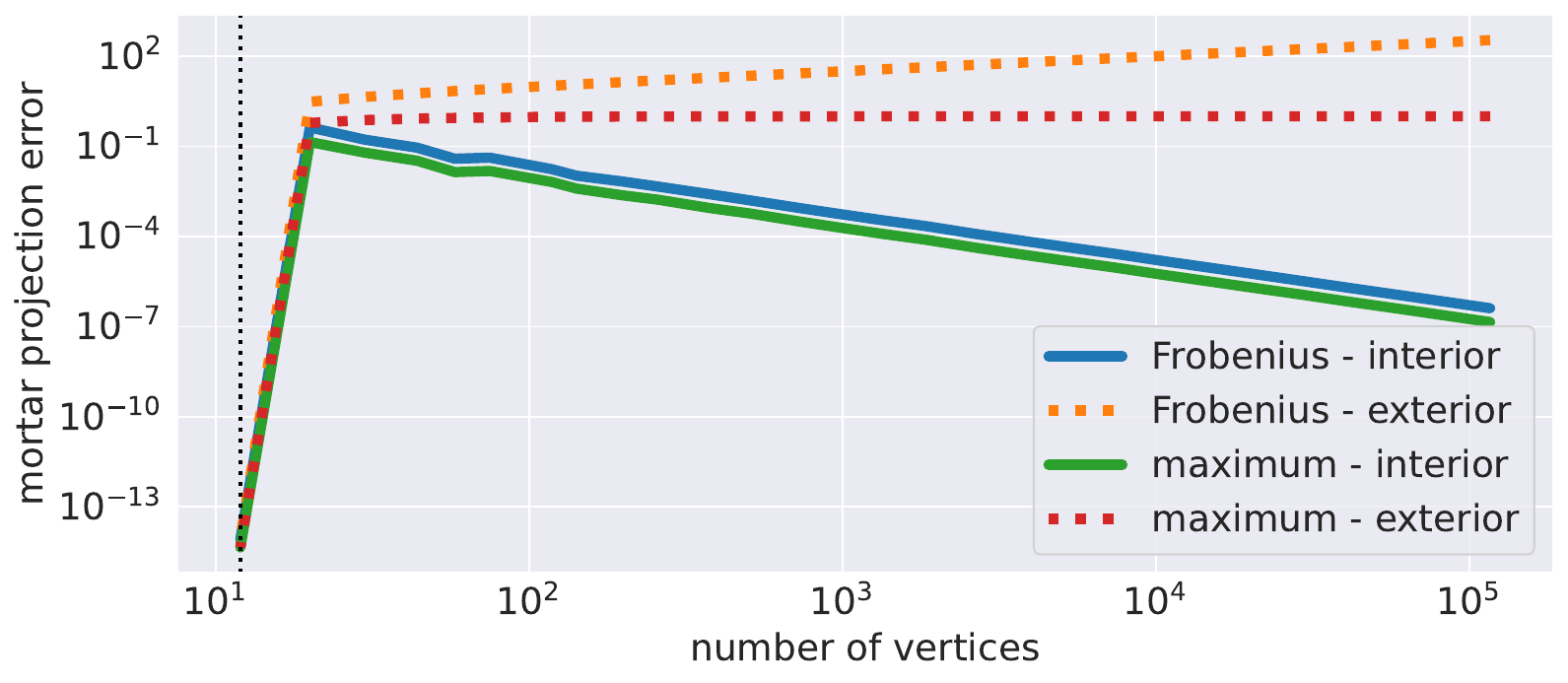}
		\caption{The projection error~\eqref{eq:error:projection} with respect to an increasingly fine `exterior' mesh. The vertical line visualises the 12 vertices of the fixed `interior' mesh.}
		\label{fig:mortarerror:meshsize}
	\end{subfigure}
	\caption{The error of the mortar projection between nonconforming meshes.}
	\label{fig:mortarerror}
\end{figure}

Figure~\ref{fig:mortarerror:perturbations} depicts the projection error, in the Frobenius and maximum norm. All four variants of the projection error decay quickly with smaller perturbations and flatten out when machine precision is reached. Hence, this experiment confirms the robustness of the PANG2 algorithm in finite precision arithmetic.

The two meshes in previous example have the same size but different configuration. Let us test the nonconforming projection with mesh refinement. Figure~\ref{fig:mortarerror:meshsize} shows the results when a fixed `interior' mesh with $h=0.5$ and 12~vertices is combined with an `exterior' mesh that has an increasingly higher resolution. The finest mesh has $h=0.0063$ and 116\,746 nodes at the unit square $[0,1] \times [0,1]$. The leftmost data point in the line plot corresponds to equal mesh sizes and yields a projection error in the order of the machine precision. When refining the `exterior' mesh, the interior projection error~\eqref{eq:error:projection:int} decreases while the exterior projection error~\eqref{eq:error:projection:ext} increases. This deterioration in projection accuracy might be due to the passage from a high-dimensional space to a low-dimensional space and back to the high-dimensional space, thus losing precision in the dimension reduction.

\subsection{Accuracy of nonconforming BEM}

It is common practice to generate meshes for the BEM based on a fixed number of elements per wavelength. While conforming BEM uses a mesh corresponding to the smallest wavelength across the interface, nonconforming BEM is applicable to meshes based on the wavenumber in each subdomain. This reduction in degrees of freedom at the same frequency is one of the key advantages of the nonconforming BEM. On the downside, nonconforming BEM introduces inaccuracies in the mortar projection. Hence, let us test the accuracy of the nonconforming BEM on a cube of size $[0,1] \times [0,1] \times [0,1]$. The incident plane wave field travels in the positive $x$-direction, i.e.,
\begin{equation}
	\uinc(\mathbf{x}) = e^{\imath \frac{2\pi f}{\cext} \mathbf{x} \cdot \mathbf{d}} \quad \text{for } \mathbf{x} \in \mathbb{R}^3 \text{ and } \mathbf{d} = \begin{bmatrix} 1 \\ 0 \\ 0 \end{bmatrix}.
	\label{eq:planewave}
\end{equation}
The physical parameters are chosen as $f=1$, $\cext=0.3$, $\cint=1.1$, $\rhoext=1$, and $\rhoint=2$. In this case, the wavelengths are $\lambdaext=0.3$ and $\lambdaint=1.1$ {\color{black} and the wavenumbers are $\kext=20.9$ and $\kext=5.7$. These values compare to the unit cube with nondimensional edge lengths equal to one}. The acoustic field is calculated on a grid of $100 \times 100$ points uniformly located in the square $[-0.5,1.5] \times [-0.5,1.5]$ on the plane $z=0.5$. No analytical solution is available for a cube and the accuracy is defined as the relative $\ell_2$ norm of the difference between the simulation and a reference solution. Specifically, the exterior PMCHWT formulation on a mesh with $\hext = \hint = 0.012$ and 49,098 vertices, corresponding to at least 25 elements per wavelength, provides the reference solution. The benchmark uses $\hext=\hint$ for conforming BEM and $\lambdaext/\hext=\lambdaint/\hint$ (a constant number of elements per wavelength) for nonconforming BEM. All meshes are independently created with different resolution.

\begin{figure}[!ht]
	\centering
	\begin{subfigure}[b]{\columnwidth}
		\includegraphics[width=\columnwidth]{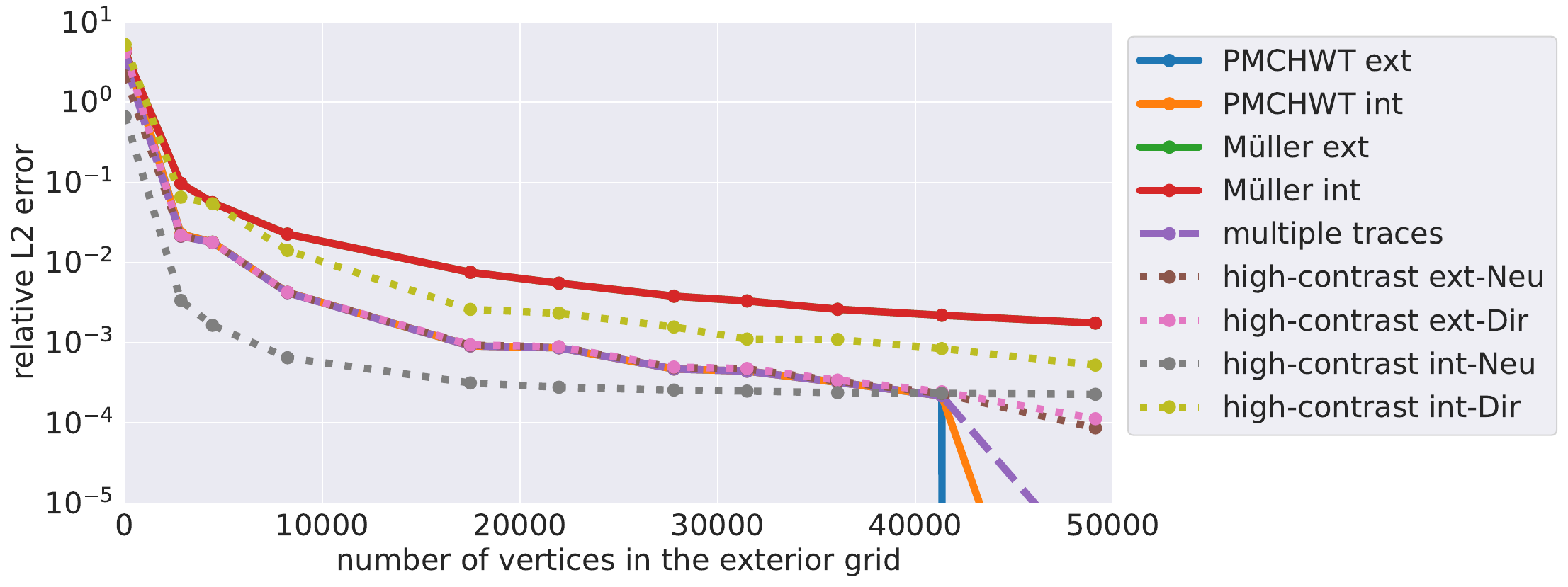}
		\caption{The accuracy of the BEM formulations for conforming meshes. The finest mesh in the convergence study corresponds to the reference simulation.}
		\label{fig:meshrefinement:conforming}
	\end{subfigure}
	\begin{subfigure}[b]{\columnwidth}
		\includegraphics[width=\columnwidth]{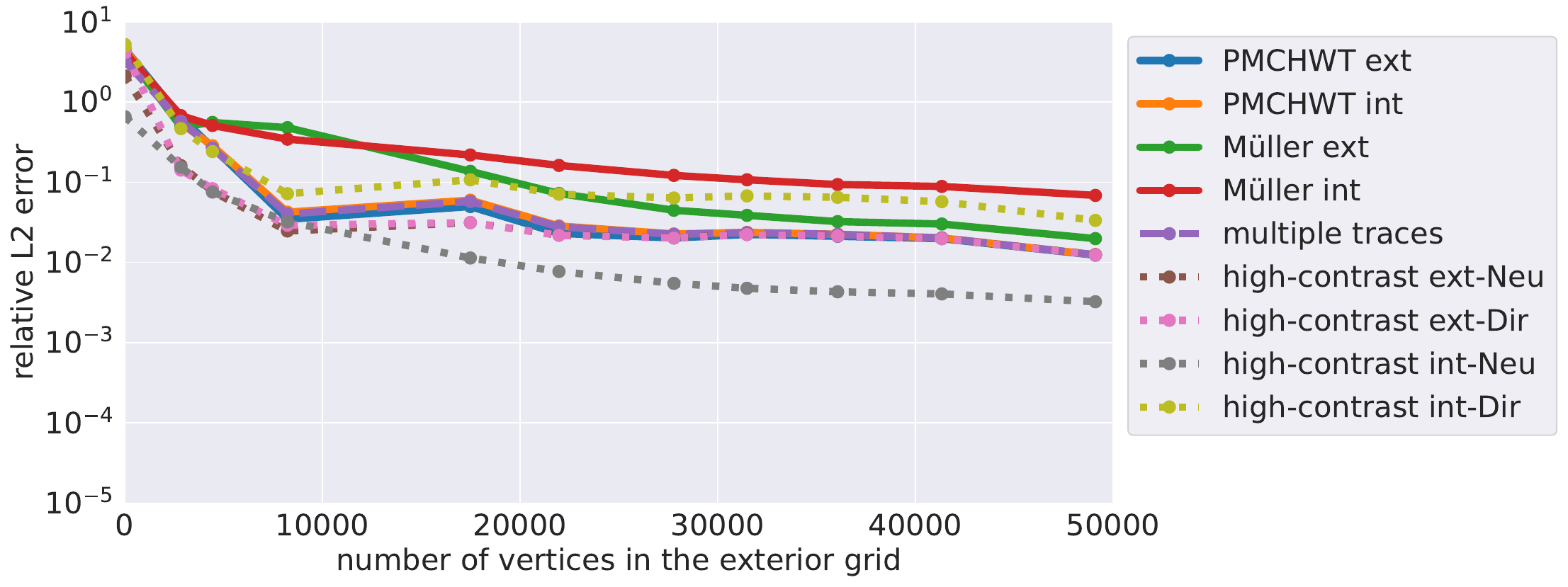}
		\caption{The accuracy of the BEM formulations for nonconforming meshes. The pair of finest meshes has 49,098 and 3788 nodes for the exterior and interior grid, respectively.}
		\label{fig:meshrefinement:nonconforming}
	\end{subfigure}
	\caption{The accuracy of the BEM with mesh refinement. The relative error is the difference with respect to the conforming exterior PMCHWT formulation with 49,098 vertices in the mesh.}
	\label{fig:meshrefinement}
\end{figure}

Figure~\ref{fig:meshrefinement} confirms that the acoustic fields improve with mesh refinement. The nonconforming BEM is not as accurate as the conforming BEM due to the additional mortar projection errors. However, it requires considerably less degrees of freedom and is thus beneficial for fast simulations that do not require high accuracy. The differences between boundary integral formulations strongly depend on the benchmark considered (cf.~\cite{wout2021benchmarking}) and its comparison is outside the scope of this study.

\FloatBarrier
\subsection{Computational efficiency}

The Bempp-cl library provides multithreaded calculations through a PyOpenCL implementation for the assembly of the discrete boundary integral operators. The PANG cannot be parallelised easily since it is an advancing front method that recursively calculates mortar contributions of neighbouring elements. However, the assembly of the mortar projection on each plane manifold of the polyhedral geometry (for example, the six faces of a cube) are independent and each task can be performed in parallel. The efficiency of the Python implementation of the mortar matrices was improved by the following approaches. Firstly, the mortar matrices satisfy $\PRextint = \PRintext^T$ so that only one of the two mortar matrices needs to be assembled. Secondly, the mortar matrix is assembled in the sparse coordinates (COO) format and then converted to a compressed-sparse-row (CSR) format for fast matrix-vector multiplications. Thirdly, the assembly of the mortar matrices on each plane manifold is implemented with the Joblib package and task-based parallelism. Fourthly, the PANG algorithm was implemented with Numba acceleration.

Let us test the efficiency of the nonconforming BEM on three benchmarks. The first is the same as in Figure~\ref{fig:meshrefinement}, the second has a three times higher frequency, and the third a three times higher material contrast, as summarised in~Table~\ref{table:efficiency:physics}. The third benchmark has a nine times lower exterior density so that the bulk modulus remains constant. While a change in density does not directly influence the mesh resolution, it does influence the conditioning of the linear system~\cite{wout2022highcontrast}. The grids consider six elements per wavelength and the mesh statistics are presented in Table~\ref{table:efficiency:mesh}.

\begin{table}[!ht]
	\caption{Physical parameters for the efficiency benchmarks {\color{black}on the unit cube}.}
	\label{table:efficiency:physics}
	\centering
	\begin{tabular}{lrrrrrrr}
		\hline\hline
		benchmark & $f$ & $\cext$ & $\cint$ & $\rhoext$ & $\rhoint$ & {\color{black}$\kext$} & {\color{black}$\kint$} \\
		\hline
		standard       & 1 & 0.3 & 1.1 & 1    & 2 & {\color{black}20.9} &  {\color{black}5.7} \\
		high frequency & 3 & 0.3 & 1.1 & 1    & 2 & {\color{black}62.8} & {\color{black}17.1} \\
		high contrast  & 1 & 0.1 & 1.1 & 0.11 & 2 & {\color{black}62.8} &  {\color{black}5.7} \\
		\hline\hline
	\end{tabular}
\end{table}

\begin{table}[!ht]
	\caption{Statistics of the meshes used in the efficiency benchmarks {\color{black}on the unit cube}.}
	\label{table:efficiency:mesh}
	\centering
	\begin{tabular}{lrrrrrrrr}
		\hline\hline
		& \multicolumn{2}{c}{both} & \multicolumn{3}{c}{conforming} & \multicolumn{3}{c}{nonconforming} \\
		benchmark & $\hext$ & $\Next$ & $\hint$ & $\Nint$ & $N_\mathrm{nodes}$ & $\hint$ & $\Nint$ & $N_\mathrm{nodes}$ \\
		\hline
		standard       & 0.05  &    2836 & 0.05  &    2836 &    5672 & 0.183 &  272 &    3108 \\
		high frequency & 0.017 & 25,302 & 0.017 & 25,302 & 50,604 & 0.061 & 2072 & 27,374 \\
		high contrast  & 0.017 & 25,302 & 0.017 & 25,302 & 50,604 & 0.183 &  272 & 25,574 \\
		\hline\hline
	\end{tabular}
\end{table}

\begin{figure}[!ht]
	\centering
	\begin{subfigure}[b]{\columnwidth}
		\includegraphics[width=.49\columnwidth]{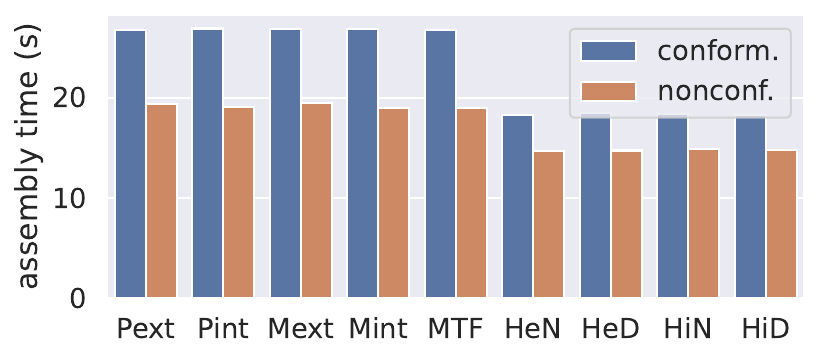}
		\includegraphics[width=.49\columnwidth]{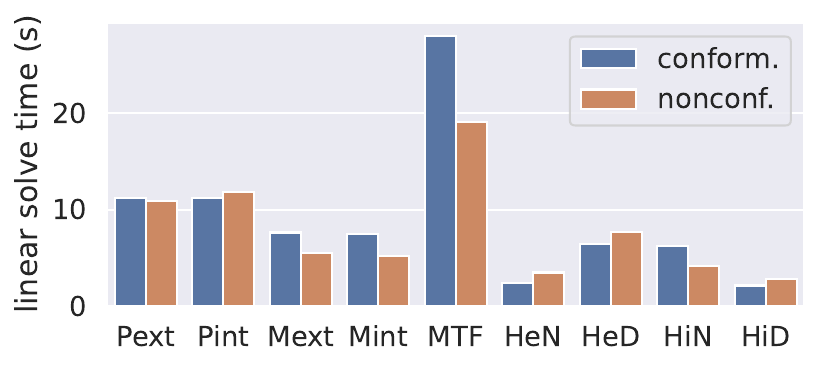}
		\includegraphics[width=.49\columnwidth]{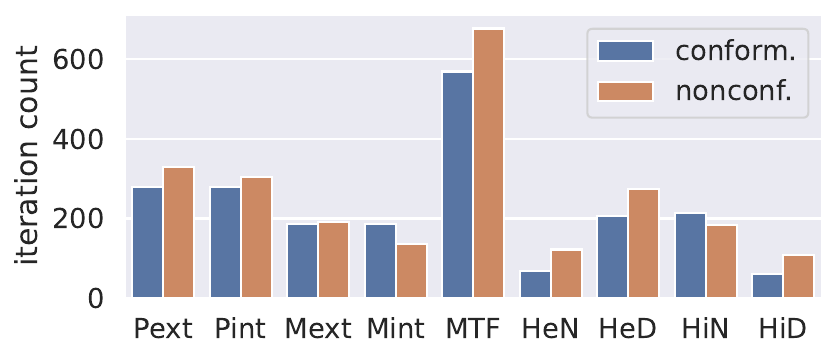}
		\includegraphics[width=.49\columnwidth]{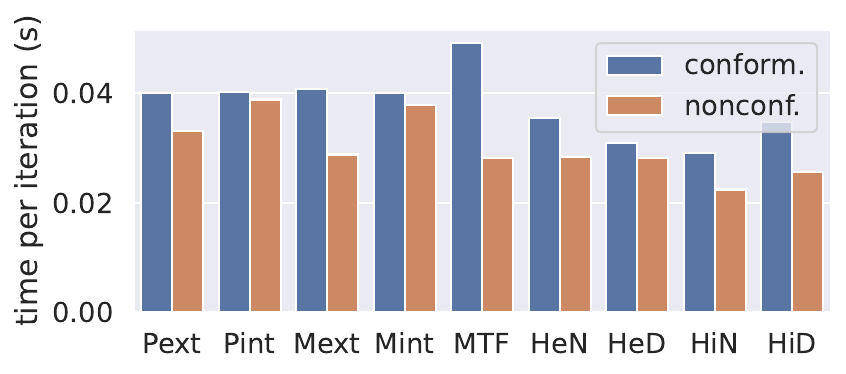}
		\caption{The benchmark with standard physical parameters.}
		\label{fig:efficiency:standard}
	\end{subfigure}
	\begin{subfigure}[b]{\columnwidth}
		\includegraphics[width=.49\columnwidth]{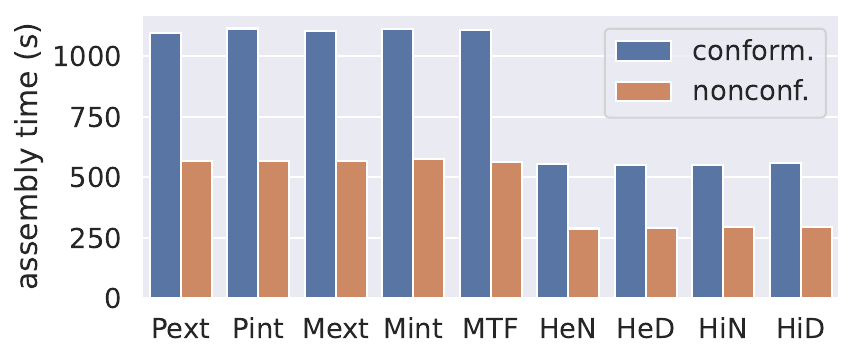}
		\includegraphics[width=.49\columnwidth]{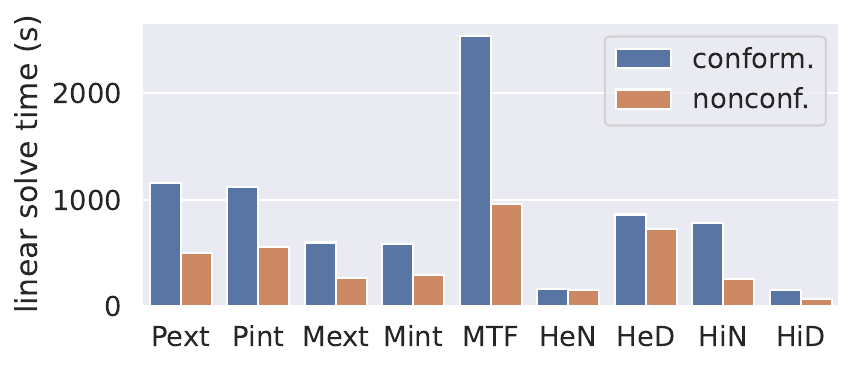}
		\includegraphics[width=.49\columnwidth]{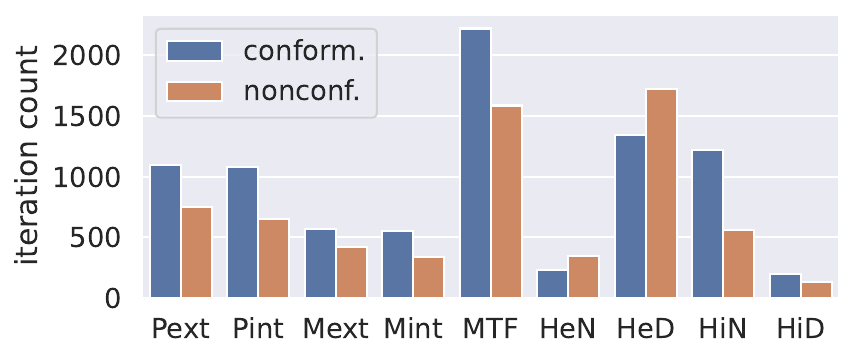}
		\includegraphics[width=.49\columnwidth]{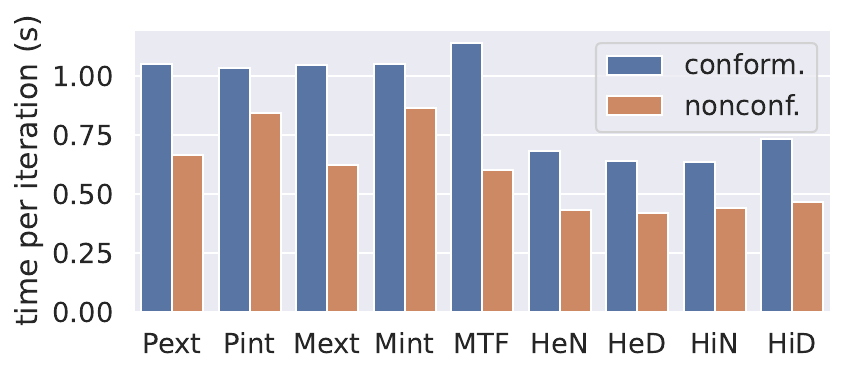}
		\caption{The benchmark with high frequency.}
		\label{fig:efficiency:highfrequency}
	\end{subfigure}
	\begin{subfigure}[b]{\columnwidth}
		\includegraphics[width=.49\columnwidth]{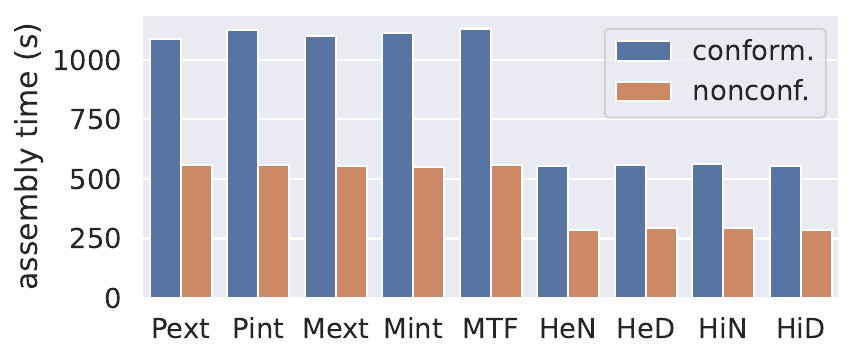}
		\includegraphics[width=.49\columnwidth]{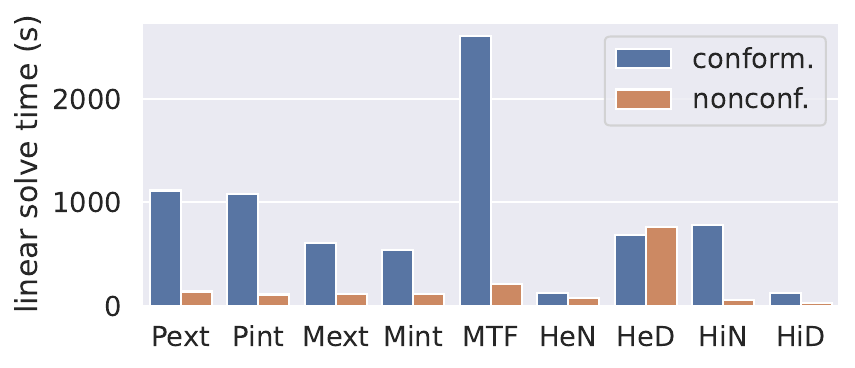}
		\includegraphics[width=.49\columnwidth]{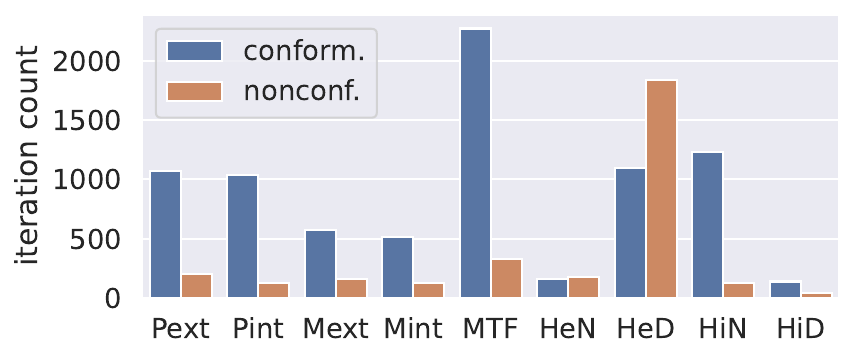}
		\includegraphics[width=.49\columnwidth]{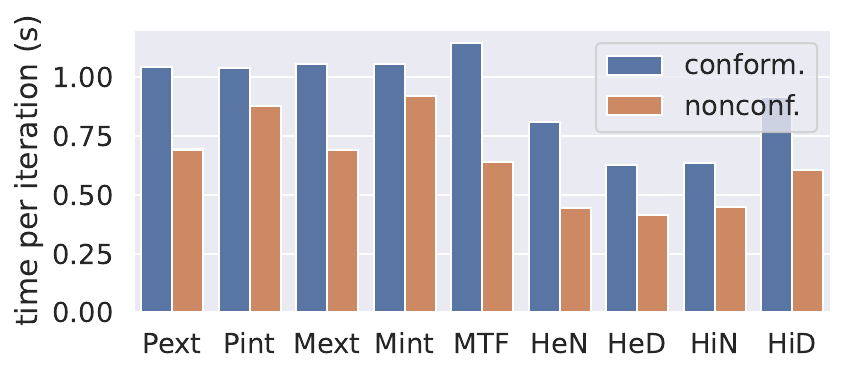}
		\caption{The benchmark with high material contrast.}
		\label{fig:efficiency:highcontrast}
	\end{subfigure}
	\caption{The efficiency statistics of the nonconforming BEM, for all formulations in Table~\ref{table:formulations}.}
	\label{fig:efficiency}
\end{figure}

Figure~\ref{fig:efficiency} presents the wall-clock time of the different benchmarks and methodologies. The direct formulations involve a full set of eight boundary integral operators, while the high-contrast formulations require four dense operators only. This halves the assembly time. The assembly time for the nonconforming BEM is almost two times faster because the interior mesh has around ten times fewer nodes, and the matrices have around a hundred times fewer elements. Hence, the assembly time of the interior operators is not significant in the overall performance. This reduction in the number of degrees of freedom also explains the differences in the time per GMRES iteration: high-contrast is faster than direct formulations and nonconforming is faster than conforming BEM. However, the differences are diluted due to each iteration's vector operations and preconditioning. The time to solve the system is dominated by the iteration count, and the differences between conforming and nonconforming BEM depend on the benchmark parameters and the boundary integral formulation. In many cases, nonconforming BEM improves the iteration count, on top of a reduction in time per iteration.

\FloatBarrier
The timings include both the boundary integral operators and the mortar matrices. Looking into the breakdown, the BEM operators dominate the overall performance. The nonconforming routines never exceed 3\% of the assembly time and 6\% of the linear solve. These small contributions are because the mortar matrix is sparse, and the boundary integral operators are dense. Notice that acceleration with fast multipole methods or hierarchical matrix compression can speed up the dense arithmetic. However, their efficiency strongly depends on the specific implementation and parameter settings. While straightforward to apply accelerators to the nonconforming BEM, they are not considered in this study to limit the parameter space of the benchmarks.

Comparing the benchmark with standard physical parameters (Fig.~\ref{fig:efficiency:standard}) to the more challenging high-frequency (Fig.~\ref{fig:efficiency:highfrequency}) and high-contrast (Fig.~\ref{fig:efficiency:highcontrast}) benchmarks, there is a significant gain in computation time for the nonconforming approach. At the challenging benchmarks, the nonconforming BEM relatively saves more nodes in the mesh, which has a quadratic influence on the size of the system matrix. Furthermore, nonconforming BEM reduces the iteration count considerably for the high-frequency and high-contrast benchmarks. Hence, the overall computation time of the nonconforming BEM is a fraction of the conforming BEM.

\subsection{FEM-BEM coupling}

Section~\ref{sec:fembem} explained the nonconforming algorithm for FEM-BEM coupling. As in the case of pure BEM, advantages include independent mesh generation and reducing the degrees of freedom for high contrasts in speed of sound across the materials. {\color{black}In contrast to pure BEM, another benefit of nonconforming FEM-BEM coupling is the opportunity to use different grid resolutions proportional to the frequency for FEM and BEM. That is, while BEM only needs six elements per wavelength~\cite{marburg2002six}, the FEM needs a higher resolution~\cite{langer2017more} to achieve accurate results in practice~\cite{aubry2022benchmark}. Furthermore, the pollution effect in the FEM at high frequencies requires disproportionally more linear elements per wavelength~\cite{babuska1997pollution}, or adaptive high-order polynomial elements~\cite{chaumont2016high}. The pollution effect has not been observed to be as restrictive for the BEM as it is for the FEM, but its understanding and quantification remains an open discussion in literature~\cite{marburg2018pollution, kreuzer2022numerical, galkowski2022helmholtz}. In any case, the nonconforming algorithm allows for generating meshes with a different number of elements per wavelength for the BEM than for the FEM, allowing for more flexibility.}

Let us test the nonconforming FEM-BEM coupling on a unit cube with a heterogeneous wavespeed. Specifically,
\begin{equation}
	\kint(x,y,z) = \kext\left( 2 + \sin(2\pi x) \sin(2\pi y) \right),
\end{equation}
$\cext=1$, $\rhoint=1$, and $\rhoext=2$. Hence, the interior wavenumber is up to three times higher than in the exterior. Let us consider an incident plane wave~\eqref{eq:planewave} with $f=2$, resulting in an exterior wavelength of $\lambdaext=0.5$ and $0.167 \leq \lambdaint \leq 0.5$ in the interior domain. {\color{black}These values are chosen proportional to the nondimensional geometry, which is again a unit cube with edge lengths equal to one.} The relatively high frequency and heterogeneous wavespeed necessitate a fine mesh for the FEM. The tetrahedral grid has 1,157,625~nodes with a maximum diameter of 0.01665 yielding at least ten elements per wavelength. The triangular surface mesh extracted from the volumetric mesh has 64,898~nodes and a maximum diameter of 0.01359, yielding at least 36~elements per exterior wavelength. The nonconforming FEM-BEM coupling uses the same FEM mesh but another BEM mesh with 5966~nodes and a maximum diameter of 0.04694, yielding at least 10~elements per exterior wavelength at the boundary.

\begin{figure}[!ht]
	\centering
	\includegraphics[width=\columnwidth]{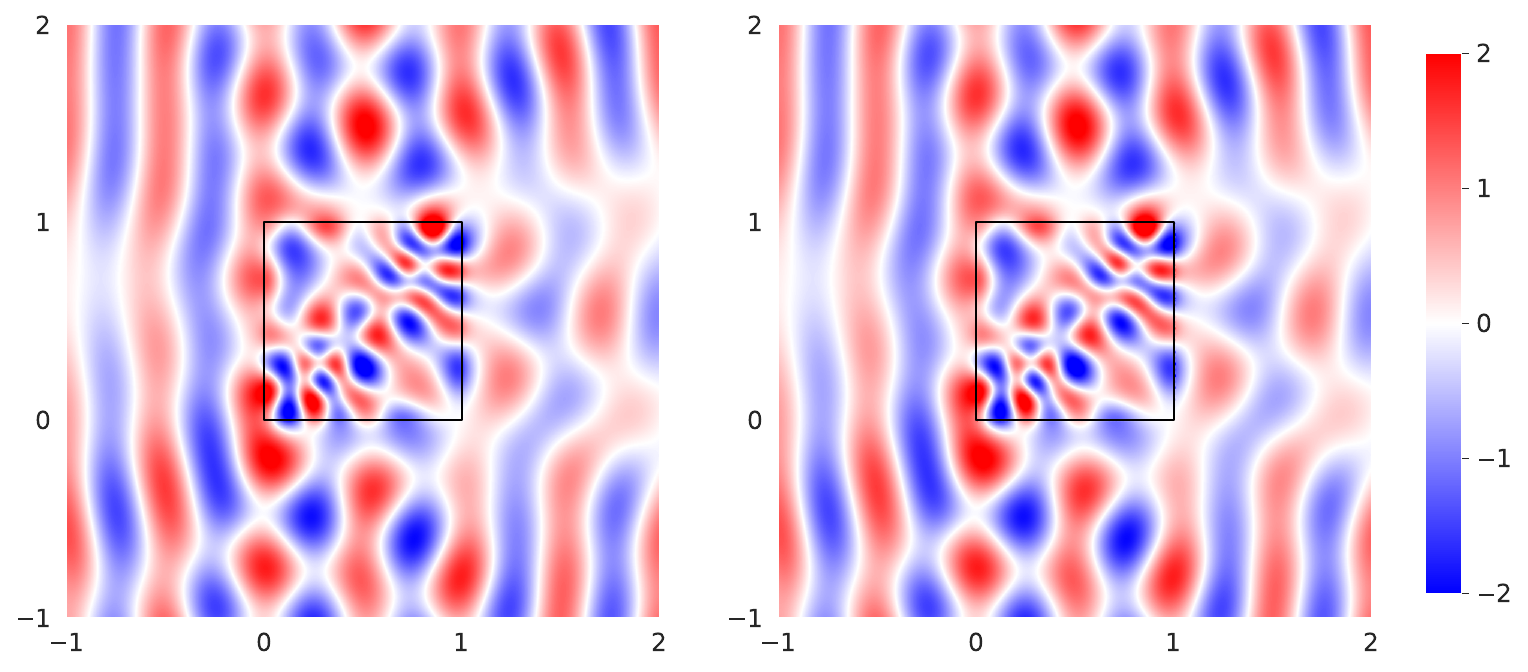}
	\caption{The real part of the acoustic field in the $(x,y)$-plane at $z=0.5$ for an incident plane wave travelling in the positive $x$-direction. The black square depicts the boundary of the unit cube. The FEM-BEM coupled system uses conforming (left) or nonconforming (right) grids.}
	\label{fig:fembem}
\end{figure}

The standard FEM-BEM coupled formulation~\eqref{eq:fembem} is discretised with P1 elements for both the BEM and FEM. The BEM part is preconditioned with an inverse mass matrix and the FEM part with an incomplete LU factorisation, implemented with the SuperLU algorithm and default parameters~\cite{li2005overview}. The linear solver converged in 12,662~iterations for the nonconforming algorithm and in 13,961~iterations for the conforming FEM-BEM coupling. The GMRES solver was restarted every thousand iterations to reduce memory consumption. The entire simulation took five hours for the nonconforming and twelve hours for the conforming FEM-BEM algorithm. The acoustic field is depicted in Figure~\ref{fig:fembem} and clearly shows the impact of a heterogeneous wavespeed in the interior. Also, the field from the conforming and nonconforming algorithm are almost identical, thus confirming the validity of the nonconforming FEM-BEM coupling.

\FloatBarrier
\subsection{Nonconforming operator preconditioning}

\begin{figure}[!ht]
	\centering
	\includegraphics[width=.9\columnwidth]{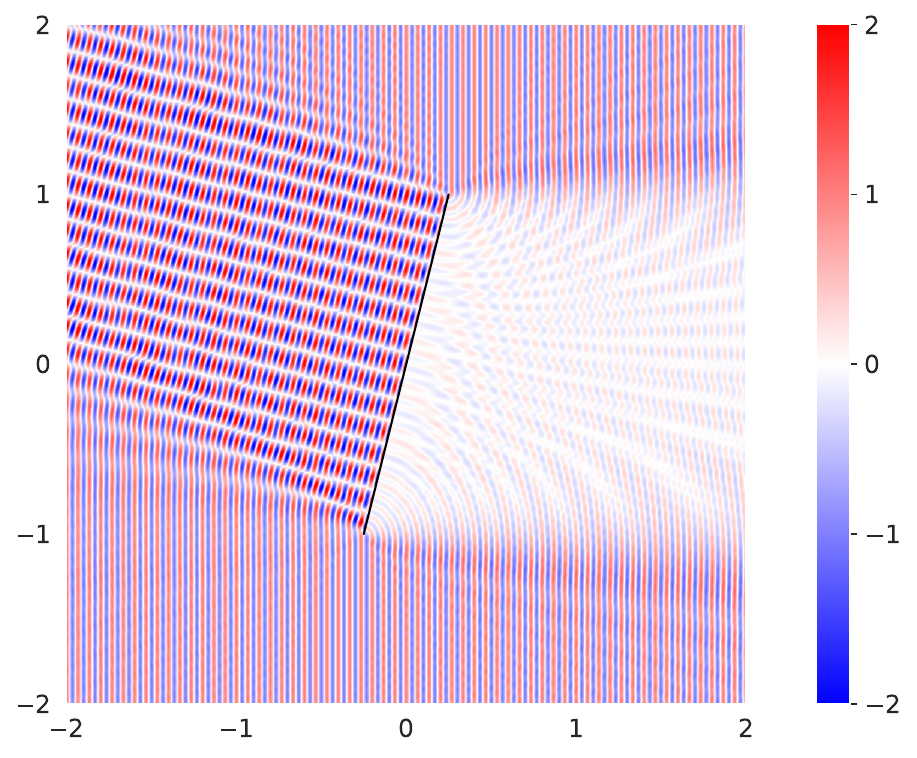}
	\caption{The real part of the acoustic field in the $(x,y)$-plane at $z=0$, calculated with nonconforming operator-preconditioned BEM. The black line depicts the slice of the screen, which is perpendicular to the visualisation plane.}
	\label{fig:screen}
\end{figure}

The nonconforming algorithm is not limited to acoustic transmission and also applies to operator preconditioning at impenetrable domains, see Section~\ref{sec:opprec}. Let us consider the Neumann problem~\eqref{eq:neumannscreen} on a rectangular screen with corners at $(-0.25, -1, -1)$, $(0.25, 1, -1)$, $(0.25, 1, 1)$, and $(-0.25, -1, 1)$; {\color{black}having nondimensional edge lengths of 2 and 2.1}. The propagating medium has $\cext=1$ and $\rhoext=1$, and the incident plane wave field propagates in the positive $x$-direction with frequency $f=15$ {\color{black}yielding $\kext=94.2$}. The surface mesh has 87,712~nodes and at least six elements per wavelength. The nonconforming operator preconditioner takes a coarse mesh with 9934~nodes and at least two elements per wavelength. Figure~\ref{fig:screen} shows the field and Table~\ref{table:screen} the timing characteristics. The single-layer operator reduces the iteration count drastically and is thus an effective preconditioner for the hypersingular operator. However, the assembly time almost doubles due to the assembly of another dense matrix. The nonconforming approach alleviates this issue: the assembly on the coarse mesh and the mortar matrices are swift, the number of iterations remains low, and the linear system is solved quickly. Nevertheless, the assembly dominates the solve time, and the overall gain in simulation time is limited compared to mass-matrix preconditioning.

\begin{table}[!ht]
	\caption{The computational efficiency of the mass-preconditioned and  opposite-order (OO) operator-preconditioned BEM on a screen: iteration count of GMRES (\#iter), wall-clock time of the simulation ($T_\mathrm{bem}$), matrix assembly ($T_\mathrm{matrix}$), linear solver ($T_\mathrm{solve}$), and each iteration ($T_\mathrm{iter}$) in seconds.}
	\label{table:screen}
	\centering
	\begin{tabular}{lrrrrrr}
		\hline\hline
		formulation & \#iter & $T_\mathrm{bem}$ & $T_\mathrm{matrix}$ & $T_\mathrm{solve}$ & $T_\mathrm{iter}$ \\
		\hline
		mass preconditioner    & 72 & 2136 & 1928 & 197 & 2.73 \\
		conforming OO prec.    &  7 & 3612 & 3534 &  67 & 9.63 \\
		nonconforming OO prec. &  7 & 2064 & 2028 &  25 & 3.53 \\
		\hline\hline
	\end{tabular}
\end{table}

\FloatBarrier
\subsection{Acoustic foam model}

As a final simulation to assess the nonconforming BEM, let us consider a large-scale benchmark on a geometry representing an acoustic foam. The object has a thickness of 2~cm for the base layer, and each pyramid has a rectangular base of size $10 \times 10$~cm and a height of 7~cm. The pyramids are replicated $16 \times 16$ times in a regular grid in the horizontal direction. The incident field
\begin{equation}
	\uinc(\mathbf{x}) = \frac{e^{\imath \frac{2\pi f}{\cext} |\mathbf{x} - \mathbf{x}_\text{source}|}}{4\pi |\mathbf{x} - \mathbf{x}_\text{source}|} \quad \text{for } \mathbf{x} \ne \mathbf{x}_\text{source} \text{ and } \mathbf{x}_\text{source} = \begin{bmatrix} 0 \\ 0 \\ 0.15 \end{bmatrix}
\end{equation}
is a point source located at 15~cm altitude in the middle of the foam, with frequency of $f=3$~kHz, which is at the higher end of the audible range. The exterior medium has material parameters $\cext=340$~m/s and $\rhoext=1.225$~kg/m$^3$ that represent air. The interior medium represents a plastic foam, where attenuation was included as a power law for the complex wavenumber $\kint = 2\pi f/\cint + \imath \alpha_\interior f/f_\alpha$ with $\alpha$ the absorption coefficient~\cite{garrett2020understanding}. The material parameters represent polyurethane~\cite{sinha2007acoustic}: $\cint=1104$~m/s, $\rhoint=1750$~kg/m$^3$, and $\alpha_\interior=7.5$~dB/cm~$=86.3$~Neper/m at $f_\alpha=2$~MHz. The fine mesh has at least 5.7~elements per wavelength in the exterior and 18.5 in the interior. The coarse mesh has at least 7~elements per wavelength in the interior domain.

\begin{figure}[!ht]
	\centering
	\includegraphics[width=\columnwidth]{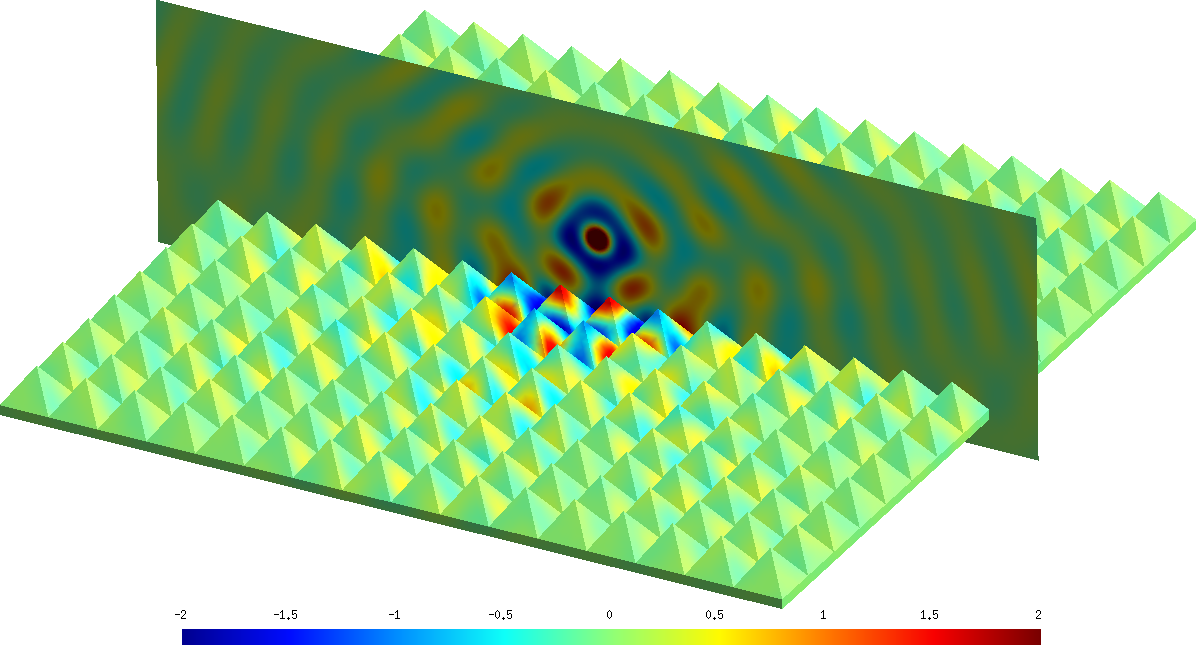}
	\caption{The real part of the acoustic field on the surface of the foam and on a slice through the point source, calculated with the nonconforming exterior PMCHWT formulation.}
	\label{fig:foam}
\end{figure}

\begin{table}[!ht]
	\caption{The number of nodes in the meshes and the BEM's efficiency on the foam model: iteration count of GMRES (\#iter), wall-clock time of the assembly ($T_\mathrm{matrix}$), linear solver ($T_\mathrm{solve}$), and each iteration ($T_\mathrm{iter}$) in seconds.}
	\label{table:foam}
	\centering
	\begin{tabular}{lrrrrrr}
		\hline\hline
		formulation & $\Next$ & $\Nint$ & \#iter & $T_\mathrm{matrix}$ & $T_\mathrm{solve}$ & $T_\mathrm{iter}$ \\
		\hline
		\multicolumn{7}{c}{conforming BEM} \\
		\hline
		PMCHWT ext & 43,469 & 43,469 & 4076 & 3478 & 17224 & 4.23 \\
		PMCHWT int & 43,469 & 43,469 & 3838 & 3453 & 15296 & 3.99 \\
		Müller ext & 43,469 & 43,469 & 4074 & 3616 & 19503 & 4.79 \\
		Müller int & 43,469 & 43,469 & 3827 & 3403 & 13973 & 3.65 \\
		multiple-traces & 43,469 & 43,469 & 8272 & 3378 & 34186 & 4.13 \\
		high-contrast ext Neu & 43,469 & 43,469 & 175 & 1696 & 299 & 1.71 \\
		high-contrast ext Dir & 43,469 & 43,469 & 2608 & 1690 & 4233 & 1.62 \\
		high-contrast int Neu & 43,469 & 43,469 & 3268 & 1743 & 5363 & 1.64 \\
		\hline
		\multicolumn{7}{c}{nonconforming BEM} \\
		\hline
		PMCHWT ext & 43,469 & 6707 & 3626 & 1745 & 6187 & 1.71 \\
		PMCHWT int & 43,469 & 6707 & 1502 & 1747 & 3017 & 2.01 \\
		Müller ext & 43,469 & 6707 & 3257 & 1783 & 5511 & 1.69 \\
		Müller int & 43,469 & 6707 & 1356 & 1751 & 2740 & 2.02 \\
		multiple-traces & 43,469 & 6707 & 6988 & 1767 & 12247 & 1.75 \\
		high-contrast ext Neu & 43,469 & 6707 & 175 & 893 & 181 & 1.04 \\
		high-contrast ext Dir & 43,469 & 6707 & 3927 & 933 & 4108 & 1.05 \\
		high-contrast int Neu & 43,469 & 6707 & 1883 & 931 & 1933 & 1.03 \\
		\hline\hline
	\end{tabular}
\end{table}

Figure~\ref{fig:foam} presents the acoustic field and Table~\ref{table:foam} the efficiency characteristics of the different boundary integral formulations. Generally speaking, the nonconforming BEM improves the computational efficiency significantly. The assembly time almost halves, each iteration is quicker, and the iteration count is often lower than conforming BEM. The high-contrast interior Dirichlet formulation did not converge to the correct solution.

Overall, the results confirm the capacity of the nonconforming algorithm to improve the efficiency of the BEM for acoustic wave transmission. For example, the simulation time of the exterior PMCHWT reduces from 5:45~hours to 2:12~hours. Also, storing the dense matrices of all eight boundary integral operators in double precision takes at least 225~GByte for the conforming BEM and 115~GByte for the nonconforming BEM.

\FloatBarrier
\section{Conclusions}

The nonconforming BEM efficiently simulates acoustic transmission with independent surface meshes at the material interface. The mortar matrices that couple nonconforming surface meshes can be calculated quickly and robustly with an advancing front algorithm. Generating the grids is modular, with a fixed number of elements per wavelength in each domain, thus drastically reducing the number of degrees of freedom. The nonconforming algorithm works for any boundary integral equation, such as the single-trace, multiple-traces, and high-contrast formulations. Furthermore, the nonconforming algorithm improves FEM-BEM coupling and operator preconditioning for models involving heterogeneous or impenetrable structures.

The salient features of the nonconforming BEM include the flexibility in mesh generation and relaxed constraints on the mesh resolution. The computational benchmarks confirm a significant improvement in computational efficiency, reducing the calculation time of matrix assembly and matrix-vector products. On top of this, many benchmarks show faster convergence of the GMRES linear solver.

The downside of the nonconforming BEM is the introduction of projection errors between the meshes, which can be a limiting factor for high-accuracy simulations. Future research will consider the more advanced techniques proposed in the nonconforming FEM literature. Furthermore, the applicability of nonconforming BEM needs to be extended to geometries that have junctions with more than two subdomains, and curved surfaces where nonmatching grids are present.

\section*{Acknowledgements}
This work was financially supported by the Agencia Nacional de Investigación y Desarrollo, Chile [FONDECYT 1230642].

\bibliographystyle{unsrt}
{\small
\bibliography{refs.bib}
}

\end{document}